\documentclass[a4paper,12pt]{article}%
\usepackage{amssymb}
\usepackage{amsfonts}
\usepackage{graphicx}
\usepackage{amsmath}%
\newtheorem{theorem}{Theorem}

\newtheorem{definition}[theorem]{Definition}

\newtheorem{lemma}[theorem]{Lemma}

\newtheorem{proposition}[theorem]{Proposition}

\newtheorem{remark}[theorem]{Remark}

\begin{document}

\title{Existence of hylomorphic solitary waves in Klein-Gordon and in
Klein-Gordon-Maxwell equations}
\author{Vieri Benci$^{\ast}$, Donato Fortunato$^{\ast\ast}$\\$^{\ast}$Dipartimento di Matematica Applicata \textquotedblleft U.
Dini\textquotedblright\\Universit\`{a} degli Studi di Pisa Universit\`{a} di Pisa\\via Filippo Buonarroti 1/c, 56127 Pisa, Italy\\e-mail: benci@dma.unipi.it\\$^{\ast\ast}$Dipartimento di Matematica \\Universit\`{a} di Bari and INFN sezione di Bari\\Via Orabona 4, 70125 Bari, Italy\\e-mail: fortunat@dm.uniba.it}
\maketitle
\tableofcontents

\textit{Dedicated to the memory of Guido Stampacchia}

\bigskip

\section{Introduction}

Roughly speaking a \textit{solitary wave} is a solution of a field equation
whose energy travels as a localized packet and which preserves this
localization in time. A solitary wave which has a non-vanishing angular
momentum is called \textit{vortex}. A \textit{soliton} is a solitary wave
which exhibits some strong form of stability so that it has a particle-like
behavior (see e.g. \cite{Ba-Be-R.}, \cite{befogranas}, \cite{raj}, \cite{vil}).\ 

To day, we know (at least) three mechanisms which might produce solitary
waves, vortices and solitons:

\begin{itemize}
\item Complete integrability, (e.g. Kortewg-de Vries equation);

\item Topological constraints, (e.g. Sine-Gordon equation);

\item Ratio energy/charge: (e.g. the nonlinear Klein-Gordon equation).
\end{itemize}

Following \cite{hylo}, the third type of solitary waves or solitons will be
called \textit{hylomorphic}. This class includes the $Q$-\textit{balls }which
are spherically symmetric solutions of the nonlinear Klein-Gordon equation
(NKG) (see \cite{Coleman86}, \cite{Ku97}) as well as solitary waves and
vortices which occur, by the same mechanism, in the nonlinear Schr\"{o}dinger
equation and in gauge theories (\cite{befo}, \cite{befo08}).

This paper is devoted to an abstract theorem which allows to prove the
existence of hylomorphic solitary waves, solitons and vortices in the (NKG)
and in the nonlinear Klein-Gordon-Maxwell equations (NKGM) .

\section{Hylomorphic solitons}

In this section we will sketch the main ideas relative to hylomorphic
solitons. They can be considered as particular states of a system modelled by
a field equation.

We assume that the state of the system is described by one or more fields
which mathematically are represented by a function
\begin{equation}
\Psi:\mathbb{R}^{N}\rightarrow V\label{lilla}%
\end{equation}
where $V$ is a finite dimensional vector space with norm $\left\vert
\ \cdot\ \right\vert _{V}$ and it is called the internal parameters space. We
will denote by $\mathcal{X}$ the set of all the states.

A state $\Psi_{0}\in\mathcal{X}$ is called solitary wave if its
evolution$\ \Psi(t)$ has the following form:%
\[
\Psi(t)=h_{t}\Psi_{0}(g_{t}x)
\]
where $h_{t}$ and $g_{t}$ are transformations on $V$ and $\mathbb{R}^{N}$
respectively and which depend continuously on $t$. A solitary wave $\Psi
_{0}\in\mathcal{X}$ is called soliton if it is \textit{orbitally stable} i.e.
if $\Psi_{0}$ $\in\Gamma$, where $\Gamma$ is a finite dimensional manifold
which is invariant and stable (see e.g. \cite{BBBM}).

In this paper, we shall consider two cases:

- Equation (NKG) (see section 4.1) where%
\[
\Psi=(\psi,\psi_{t})\in\mathbb{C}^{2}.
\]
- Equation (NKGM) (see section 6.1) where%
\[
\Psi=(\psi,\psi_{t},\phi,\phi_{t},\mathbf{A},\mathbf{A}_{t})\in\mathbb{C}%
^{2}\times\mathbb{R}^{8}.
\]

The existence and the properties of hylomorphic solitons are guaranteed by the
interplay between \textit{energy }$E$ and another integral of motion which, in
the general case, is called \textit{hylenic charge} and it will be denoted
by\textit{\ }$H.$

Thus, the most general equations for which it is possible to have hylomorphic
solitons need to have the following features:

\begin{itemize}
\item \textbf{A-1. }\textit{The equations are variational namely they are the
Euler-Lagrange equations relative to a Lagrangian density $\mathcal{L}$%
}$\left[  \Psi\right]  $\textit{.}

\item \textbf{A-2. }\textit{The equations are invariant for time and space
translations, namely $\mathcal{L}$ does not depend explicitly on }$t$
\textit{and} $x. $

\item \textbf{A-3. }\textit{The equations are invariant for a }$S^{1}%
$\textit{action, namely $\mathcal{L}$ does not depend explicitly on the phase
of the field }$\Psi$ \textit{which is supposed to be complex valued (or at
lest to have some complex valued component).}
\end{itemize}

More exactly, in (NKG), we have the following $S^{1}$ action
\[
T_{\theta}\Psi=T_{\theta}(\psi,\psi_{t})=(e^{i\theta}\psi,e^{i\theta}\psi
_{t}),\ \ \theta\in\mathbb{R}/\left(  2\pi\mathbb{Z}\right)  =S^{1}%
\]
and in (NKGM) we have%
\[
T_{\theta}\Psi=T_{\theta}(\psi,\psi_{t},\phi,\phi_{t},\mathbf{A}%
,\mathbf{A}_{t})=(e^{i\theta}\psi,e^{i\theta}\psi_{t},\phi,\phi_{t}%
,\mathbf{A},\mathbf{A}_{t}).
\]
Solitary waves or solitons for equations satisfying \textbf{A-1} and
\textbf{A-2 }and having null momentum are called \textit{stationary waves }or
\textit{stationary solitons}.

By Noether theorem assumptions \textbf{A-1} and \textbf{A-2} guarantee the
conservation of the energy $E\left(  \Psi\right)  $ and of the momentum
$\mathbf{P}\left(  \Psi\right)  $ (see e.g. \cite{befogranas}), while
\textbf{A-1} and \textbf{A-3 }guarantee the conservation of another integral
of motion which we call \textit{hylenic charge} $H\left(  \Psi\right)  $ (see
\cite{hylo}).

The quantity
\begin{equation}
\Lambda\left(  \Psi\right)  =\frac{E\left(  \Psi\right)  }{\left\vert H\left(
\Psi\right)  \right\vert },\label{lambda}%
\end{equation}

\noindent which is an invariant of the motion having the dimension of energy,
is called \textit{hylomorphy ratio.}

We now set%

\begin{equation}
m=\;\underset{\varepsilon\rightarrow0}{\lim}\;\underset{\Psi\in\mathcal{X}%
_{\varepsilon}}{\inf}\frac{E\left(  \Psi\right)  }{\left\vert H\left(
\Psi\right)  \right\vert }\label{brutta}%
\end{equation}
where
\begin{equation}
\mathcal{X}_{\varepsilon}=\left\{  \Psi\in\mathcal{X}:\forall x,\ \left\Vert
\Psi(x)\right\Vert _{V}<\varepsilon\right\}  .\label{icsep}%
\end{equation}

Now let $\Psi(t)$ be the evolution of a state such that $\Lambda\left(
\Psi(0)\right)  =\lambda<m;$ then, $\Lambda\left(  \Psi(t)\right)  =\lambda
\ $for all $t,$ and, by definition of $m,$ we have that%
\[
\underset{t\rightarrow\infty}{\lim\inf}\left\Vert \Psi(t)\right\Vert _{V}>0.
\]
Thus it it possible that $\Psi(t)$ tends to a nontrivial stable configuration.

Now let $\sigma$ be a real number and $\Psi$ be a state such that%
\begin{equation}
H\left(  \Psi\right)  =\sigma\text{ and }E\left(  \Psi\right)  =\min\left\{
E\left(  v\right)  :H(v)=\sigma\right\} \label{set}%
\end{equation}

and denote by $\Gamma_{\sigma}$ the set of such minimizers $\Psi,$ namely%
\[
\Gamma_{\sigma}=\left\{  \Psi:\Psi\text{ satisfies (\ref{set})}\right\}  .
\]

Observe that by \textbf{A-2 }the energy is a constant of the motion, then
$\Gamma_{\sigma}$ is an invariant set.

Now we give the following definition

\begin{definition}
\label{hys}A stationary wave $\Psi_{0}$ is called hylomorphic wave if
\begin{equation}
\Psi_{0}\in\Gamma_{\sigma}\text{ for some }\sigma.\label{mm}%
\end{equation}
Moreover $\Psi_{0}$ is called hylomorphic soliton if it satisfies (\ref{mm})
and if $\Gamma_{\sigma}$ is a manifold with
\[
dim(\Gamma_{\sigma})<\infty\text{ and }\Gamma_{\sigma}\text{ is stable}%
\]

\end{definition}

\begin{remark}
In the examples considered in this paper, the Lagrangian \textit{$\mathcal{L}
$}$\left[  \Psi\right]  $ is invariant for an action of the Poincar\'{e}
group. In particular, if the Lagrangian is invariant for the action of a
Lorentz boost, then the existence of stationary waves and stationary solitons
implies the existence of travelling (with velocity $\mathbf{v}$, $\left\vert
\mathbf{v}\right\vert <c)$ waves and travelling solitons respectively (see
e.g. \cite{befogranas}).
\end{remark}

\section{An abstract theorem}

In many situations the energy $E$ and the charge $H$ have the following form%
\begin{align}
E(u,\omega)  & =J(u)+\omega^{2}K(u),\label{forma}\\
H(u,\omega)  & =2\omega K(u).\label{formina}%
\end{align}
where $\omega\in\mathbb{R}$ and $J$ and $K$ are as follows:%
\begin{align*}
J(u)  & =\frac{1}{2}\left\langle L_{1}u,u\right\rangle +N_{1}(u)\\
K(u)  & =\frac{1}{2}\left\langle L_{0}u,u\right\rangle +N_{0}(u)
\end{align*}
where $L_{i}:X\rightarrow X^{\prime}$ ($i=0,1)$ are linear continuous
operators and $N_{i}$ ($i=0,1)$ are differentiable functionals defined on a
Hilbert space $X$ with a norm equivalent to the following one%
\[
\left\Vert u\right\Vert ^{2}=\left\langle L_{1}u,u\right\rangle .
\]
Here $\left\langle ,\right\rangle $ denotes the duality between $X$ and
$X^{\prime}.$

The existence of solitary waves for the field equations we are interested in
lead to study the following abstract eigenvalue problem:
\begin{equation}
J^{\prime}(u)=\omega^{2}K^{\prime}(u).\label{JK}%
\end{equation}
where $J^{\prime}$ and $K^{\prime}$ denote the differentials of $J$ and $K. $

The most natural way to solve this problem consists in minimizing $J(u)$ on
the manifold $\left\{  u:K(u)=const.\right\}  $. However the assumptions which
allow such a minima to exist are not adequate for the problems which we want
to consider. For this reason we adopt a different variational principle, which
permits also to get the existence of particular solitary waves, namely of
hylomorphic waves (see Definition \ref{hys}).

We set for $\sigma>0$
\[
M_{\sigma}=\left\{  \left(  u,\omega\right)  \in X\times\mathbb{R}%
^{+}:H(u,\omega)=\sigma\right\}  .
\]
The variational principle is contained in the following simple result:

\begin{theorem}
\label{VP}The critical points $\left(  u,\omega\right)  $ of $E$ on
$M_{\sigma}$ solve the problem (\ref{JK}).
\end{theorem}

\textbf{Proof.} Let $\left(  u,\omega\right)  \in M_{\sigma}$ be a critical
point of $E$ on $M_{\sigma}.$ Then there exists $\lambda$ real such that%
\[
\left\{
\begin{array}
[c]{c}%
\partial_{u}E\left(  u,\omega\right)  =\lambda\partial_{u}H\left(
u,\omega\right) \\
\partial_{\omega}E\left(  u,\omega\right)  =\lambda\partial_{\omega}H\left(
u,\omega\right)
\end{array}
\right.
\]
These equations can be written more explicitly%
\[
\left\{
\begin{array}
[c]{c}%
J^{\prime}\left(  u\right)  +\omega^{2}K^{\prime}\left(  u\right)
=\lambda\omega K^{\prime}\left(  u\right) \\
2\omega K\left(  u\right)  =\lambda K\left(  u\right)
\end{array}
\right.
\]
From the second equation we have $\lambda=2\omega$ and substituting in the
first one, we get that $\left(  u,\omega\right)  $ solves problem (\ref{JK}).

$\square$

\ 

The utility of Theorem \ref{VP} relies on the fact that the existence of
critical points of $E$ on $M_{\sigma}$ is guaranteed by an assumption (see
assumption (\ref{HH})), which in many physical problems is the natural one.
Moreover in some cases this assumption guarantees the stability of the solutions.

We make the following assumptions:

\begin{itemize}
\item (H1)$\ $ $J$ $\geq0$ and $J$ is coercive on $M_{\sigma}$, namely for any
sequence $\left(  u_{n},\omega_{n}\right)  \in M_{\sigma}$ we have that
($J(u_{n})$ bounded)$\ \Rightarrow$ ($u_{n}$ bounded) .

\item (H2) The differentials $N_{0}^{\prime},N_{1}^{\prime}$ of $N_{0},N_{1}$
satisfy the following compactness properties:
\end{itemize}

$N_{0}^{\prime}:X\rightarrow X^{\prime}$ is compact. Moreover, if $u_{n}$
converges weakly in $X,$ then
\begin{equation}
\left\langle N_{1}^{\prime}(u_{n})-N_{1}^{\prime}(u_{m}),u_{n}-u_{m}%
\right\rangle \rightarrow0\text{ as }n,m\rightarrow\infty\label{stare}%
\end{equation}

\begin{itemize}
\item (H3) $K(u)\geq0\ $for$\ $all $u\ $and $K(u)\neq0$ for some $u\in X.$
\end{itemize}

We shall prove the following theorem:

\begin{theorem}
\label{abstract}Assume (H1,2,3) and that there is $\bar{u},$ such that%
\begin{equation}
0<\frac{J(\bar{u})}{K(\bar{u})}<m^{2}\label{HH}%
\end{equation}
where%
\begin{equation}
m^{2}=\inf\frac{\left\langle L_{1}u,u\right\rangle }{\left\langle
L_{0}u,u\right\rangle }>0.\label{norma}%
\end{equation}
Then there exists a non empty, open set $\Sigma\subset\mathbb{R}$ such that,
for any $\sigma\in\Sigma,$ $E$ has a minimizer $\left(  u_{0},\omega
_{0}\right)  $ on $M_{\sigma}$ with $0<\omega_{0}^{2}<m^{2}.$
\end{theorem}

As an immediate consequence of Theorem \ref{VP} and Theorem \ref{abstract} we get

\begin{theorem}
\label{final}Under the assumptions of Theorem \ref{abstract} there exists a
non empty, open set $\Sigma\subset\mathbb{R}$ such that, for any $\sigma
\in\Sigma$ problem (\ref{JK}) has a solution $\left(  u,\omega\right)  ,$ such
that $0<\omega^{2}<m^{2}$, $H(u,\omega)=\sigma$ and which is a minimizer of
$E$ on $M_{\sigma}.$
\end{theorem}

We set, for $\omega>0$ and $K\left(  u\right)  >0,$
\[
\Lambda\left(  u,\omega\right)  =\frac{E\left(  u,\omega\right)  }{H\left(
u,\omega\right)  }=\frac{1}{2}\left(  \frac{J\left(  u\right)  }{K\left(
u\right)  }\cdot\frac{1}{\omega}+\omega\right)  .
\]

\begin{remark}
In this paper we will apply theorem \ref{final} in three cases. In these
cases, $E$ and $H$ will represent respectively the energy and the hylenic
charge, $\Lambda$ is the hylomorphy ratio and $m$ in (\ref{norma}) coincides
with the constant defined by (\ref{brutta}).
\end{remark}

In order to prove Theorem \ref{abstract}, we need several lemmas.

\begin{lemma}
\label{arabel}If $J$ , $K\geq0,$ then the following assertions are equivalent:

\begin{itemize}
\item (a) there is $\bar{u}\in X,$ such that%
\begin{equation}
0<\frac{J(\bar{u})}{K(\bar{u})}<m^{2}.
\end{equation}

\item (b) there exist$\ \bar{u}\in X,$ $\bar{\omega}>0$ such that%
\begin{equation}
\Lambda\left(  \bar{u},\bar{\omega}\right)  <m.\label{arabella}%
\end{equation}

\end{itemize}
\end{lemma}

\textbf{Proof.} (a)$\Rightarrow$(b) If we take $\bar{\omega}=\sqrt
{\frac{J(\bar{u})}{K(\bar{u})}},$ we have that%
\[
\Lambda\left(  \bar{u},\bar{\omega}\right)  =\frac{1}{2}\left(  \frac
{J(\bar{u})}{K(\bar{u})}\cdot\frac{1}{\bar{\omega}}+\bar{\omega}\right)
=\sqrt{\frac{J(\bar{u})}{K(\bar{u})}}<m.
\]

(b)$\Rightarrow$(a) If $\frac{1}{2}\left(  \frac{J(\bar{u})}{K(\bar{u})}%
\cdot\frac{1}{\bar{\omega}}+\bar{\omega}\right)  <m,$ then%
\[
\frac{J(\bar{u})}{K(\bar{u})}<2m\bar{\omega}-\bar{\omega}^{2}\leq
\ \underset{\omega\geq0}{\max}\left(  2m\omega-\omega^{2}\right)  =m^{2}.
\]

$\square$

\begin{lemma}
\label{zero}Assume $J,K\geq0$ and let $\left(  u_{n},\omega_{n}\right)  $ be a
sequence in $M_{\sigma},$ $\sigma>0,$ with $\Lambda\left(  u_{n},\omega
_{n}\right)  $ bounded. Then the sequences $\omega_{n}$and $J(u_{n})$ are bounded.
\end{lemma}

The proof is trivial.

We now set%
\[
\hat{c}=\underset{\omega\geq m,u\in X}{\inf}\ \Lambda\left(  u,\omega\right)
.
\]

\begin{lemma}
\label{artemide}Assume that $J,K\geq0$ and let $\left(  u_{n},\omega
_{n}\right)  $ be a sequence in $M_{\sigma},$ $\sigma>0,$ such that
$\Lambda\left(  u_{n},\omega_{n}\right)  \rightarrow c<\hat{c}.$ Then (up to a
subsequence).
\[
\omega_{n}\rightarrow\omega_{0}<m.
\]

\end{lemma}

\textbf{Proof. }Let $\left(  u_{n},\omega_{n}\right)  $ be a sequence in
$M_{\sigma},\sigma>0,$ such that $\Lambda\left(  u_{n},\omega_{n}\right)
\rightarrow c<\hat{c}.$ Since $\Lambda\left(  u_{n},\omega_{n}\right)  $ is
bounded, by Lemma \ref{zero}, $\omega_{n}$ is bounded and hence, up to a
subsequence, $\omega_{n}\rightarrow\omega_{0}.$ We have to prove that
$\omega_{0}<m.$ We argue indirectly and assume that $\omega_{n}=m_{1}%
+\delta_{n}$ with $\delta_{n}\rightarrow0$ and $m_{1}\geq m.$ Since
$\omega_{n}$ and $\Lambda\left(  u_{n},\omega_{n}\right)  $ are bounded$,$
also $\frac{J\left(  u_{n}\right)  }{K\left(  u_{n}\right)  }$ is bounded,
then easy calculations give
\begin{align*}
\Lambda\left(  u_{n},m_{1}+\delta_{n}\right)   & =\frac{1}{2}\left(
\frac{J\left(  u_{n}\right)  }{K\left(  u_{n}\right)  }\cdot\frac{1}%
{m_{1}+\delta_{n}}+m_{1}+\delta_{n}\right)  =\\
& \frac{1}{2}\left(  \frac{J\left(  u_{n}\right)  }{m_{1}K\left(
u_{n}\right)  }\left(  1+\frac{\delta_{n}}{m_{1}}\right)  ^{-1}+m_{1}%
+\delta_{n}\right) \\
& =\Lambda\left(  u_{n},m_{1}\right)  +O\left(  \delta_{n}\right)  .
\end{align*}
Then%
\begin{align*}
c  & =\ \underset{n\rightarrow\infty}{\lim}\Lambda\left(  u_{n},\omega
_{n}\right)  =\underset{n\rightarrow\infty}{\lim}\Lambda\left(  u_{n}%
,m_{1}+\delta_{n}\right)  =\underset{n\rightarrow\infty}{\lim}\left(
\Lambda\left(  u_{n},m_{1}\right)  +O\left(  \delta_{n}\right)  \right) \\
& \geq\underset{\omega\geq m,u\in X}{\inf}\ \Lambda\left(  u,\omega\right)
=\hat{c},
\end{align*}
contradicting our assumption.

$\square$

\begin{lemma}
\label{cinzia}Assume (H1,2,3). Then for any $\sigma>0,$ $\Lambda$ satisfies
$PS$ in $M_{\sigma}\ $under the level $\hat{c},$ namely, if $\left(
u_{n},\omega_{n}\right)  $ is a sequence in $M_{\sigma}$ such that
\begin{align}
\Lambda\left(  u_{n},\omega_{n}\right)   & \rightarrow c<\hat{c}\label{p}\\
\left.  d\Lambda\right\vert _{M_{\sigma}}\left(  u_{n},\omega_{n}\right)   &
\rightarrow0\ ,\ \label{pp}%
\end{align}
then $\left(  u_{n},\omega_{n}\right)  $ has a converging subsequence.
\end{lemma}

Proof. Let $\left(  u_{n},\omega_{n}\right)  $ be a sequence in $M_{\sigma}$
satisfying (\ref{p}) and (\ref{pp}). By Lemma \ref{zero} $J\left(
u_{n}\right)  $ is bounded. Then, by the coercivity of $J$ on $M_{\sigma}$, we
deduce that $u_{n}$ weakly converges (up to a subsequence) to $u_{0}\in$ $X$.
Using Lemma \ref{artemide}, up to a subsequence, we get that
\begin{equation}
\omega_{n}\rightarrow\omega_{0}<m.\label{bo}%
\end{equation}
Now we prove that $u_{n}$ converges strongly to $u_{0}.$

By (\ref{pp}) we have that there exists a sequence of real numbers
$\lambda_{n}$ such that
\[
\left\{
\begin{array}
[c]{c}%
\partial_{u}E\left(  u_{n},\omega_{n}\right)  =\lambda_{n}\partial_{u}H\left(
u_{n},\omega_{n}\right)  +\varepsilon_{n}\\
\partial_{\omega}E\left(  u_{n},\omega_{n}\right)  =\lambda_{n}\partial
_{\omega}H\left(  u_{n},\omega_{n}\right)  +\eta_{n}%
\end{array}
\right.
\]

where $\varepsilon_{n}\rightarrow0$ in $X^{\prime}$ and $\eta_{n}\rightarrow0$
in $\mathbb{R}.$ These equations can be written more explicitely as follows:%
\begin{equation}
\left\{
\begin{array}
[c]{c}%
J^{\prime}\left(  u_{n}\right)  +\omega_{n}^{2}K^{\prime}\left(  u_{n}\right)
=\lambda_{n}\omega_{n}K^{\prime}\left(  u_{n}\right)  +\varepsilon_{n}\\
2\omega_{n}K\left(  u_{n}\right)  =\lambda_{n}K\left(  u_{n}\right)  +\eta_{n}%
\end{array}
.\right. \label{uno}%
\end{equation}
By the second equation we get
\[
\lambda_{n}=2\omega_{n}-\frac{\eta_{n}}{K\left(  u_{n}\right)  }=2\omega
_{n}-\frac{2\eta_{n}\omega_{n}}{\sigma};
\]
replacing $\lambda_{n}$ in the first equation, we get%
\[
J^{\prime}\left(  u_{n}\right)  -\omega_{n}^{2}K^{\prime}\left(  u_{n}\right)
=-\frac{2\eta_{n}\omega_{n}^{2}}{\sigma}K^{\prime}\left(  u_{n}\right)
+\varepsilon_{n}.
\]
This equation can be rewritten as follows%
\begin{equation}
L_{1}u_{n}-\omega_{0}^{2}L_{0}u_{n}=-N_{1}^{\prime}(u_{n})+\omega_{n}^{2}%
N_{0}^{\prime}\left(  u_{n}\right)  +\delta_{n}\label{clarabella}%
\end{equation}
where%
\[
\delta_{n}=-(\omega_{0}^{2}-\omega_{n}^{2})L_{0}u_{n}-\frac{2\eta_{n}%
\omega_{n}^{2}}{\sigma}K^{\prime}\left(  u_{n}\right)  +\varepsilon_{n}.
\]
Since $u_{n}\ $is bounded,$\ L_{0}u_{n}$ and $K^{\prime}\left(  u_{n}\right)
$ are bounded;$\ $then $\delta_{n}\rightarrow0.$

Replacing in (\ref{clarabella}) $n$ with $m$%
\begin{equation}
L_{1}(u_{m})-\omega_{0}^{2}L_{0}(u_{m})=-N^{\prime}(u_{m})+\omega_{m}^{2}%
N_{0}^{\prime}\left(  u_{m}\right)  +\delta_{m}\label{rebella}%
\end{equation}
and, subtracting (\ref{rebella}) from (\ref{clarabella}), we get%
\[
L_{1}(u_{n}-u_{m})-\omega_{0}^{2}L_{0}(u_{n}-u_{m})=
\]%
\begin{equation}
=N_{1}^{\prime}(u_{m})-N_{1}^{\prime}(u_{n})+\omega_{n}^{2}N_{0}^{\prime
}\left(  u_{n}\right)  -\omega_{m}^{2}N_{0}^{\prime}\left(  u_{m}\right)
+\delta_{n}-\delta_{m}.\label{bek}%
\end{equation}

By (H2) and since $u_{n}$ is bounded, we easily get%
\begin{equation}
\left\langle N_{1}^{\prime}(u_{m})-N_{1}^{\prime}(u_{n})+\omega_{n}^{2}%
N_{0}^{\prime}\left(  u_{n}\right)  -\omega_{m}^{2}N_{0}^{\prime}\left(
u_{m}\right)  ,u_{n}-u_{m}\right\rangle \rightarrow0.\label{bat}%
\end{equation}

By (\ref{norma}) we have that
\begin{align}
& \left\langle L_{1}\left(  u_{n}-u_{m}\right)  ,u_{n}-u_{m}\right\rangle
-\omega_{0}^{2}\left\langle L_{0}\left(  u_{n}-u_{m}\right)  ,u_{n}%
-u_{m}\right\rangle \nonumber\\
& \geq\left\langle L_{1}\left(  u_{n}-u_{m}\right)  ,u_{n}-u_{m}\right\rangle
-\frac{\omega_{0}^{2}}{m^{2}}\left\langle L_{1}\left(  u_{n}-u_{m}\right)
,u_{n}-u_{m}\right\rangle \nonumber\\
& \geq\left(  1-\frac{\omega_{0}^{2}}{m^{2}}\right)  \left\Vert u_{n}%
-u_{m}\right\Vert ^{2}.\label{bit}%
\end{align}

Thus, multiplying both sides of (\ref{bek}) by $u_{n}-u_{m}$ and using
(\ref{bat}), (\ref{bit}), we get%
\begin{equation}
\varepsilon_{n,m}\geq\left(  1-\frac{\omega_{0}^{2}}{m^{2}}\right)  \left\Vert
u_{n}-u_{m}\right\Vert ^{2}\text{ where }\varepsilon_{n,m}\rightarrow
0.\label{cauchy}%
\end{equation}

Since $\omega_{0}<m$ (see (\ref{bo})), by (\ref{cauchy}) $u_{n}$ is a Cauchy
sequence in $X$ $.$

$\square$

\begin{lemma}
\label{stupido}If assertion (a) (or (b)) in lemma \ref{arabel} holds, then
\[
\hat{c}<m.
\]

\end{lemma}

\textbf{Proof.} By (a) in lemma \ref{arabel} we have, for a suitable $\bar
{u}\in X,$ $\frac{J(\bar{u})}{K(\bar{u})}<m^{2}$.\ Then, by definition of
$\hat{c},$
\[
\hat{c}\leq\Lambda\left(  \bar{u},m\right)  =\frac{1}{2}\left(  \frac
{J(\bar{u})}{K(\bar{u})}\cdot\frac{1}{m}+m\right)  <m
\]

$\square$

Now we set%
\begin{equation}
\Sigma=\left\{  \sigma>0:\underset{\left(  u,\omega\right)  \in M_{\sigma}%
}{\inf}\Lambda\left(  u,\omega\right)  <\hat{c}\right\}  .\label{defined}%
\end{equation}
\bigskip

\bigskip The following Lemma guarantees that the set $\Sigma$ is not empty.

\begin{lemma}
\label{pomponia} If assertion (a) (or (b)) in lemma \ref{arabel} holds, then
\[
\underset{\left(  u,\omega\right)  \in X\times\mathbb{R}^{+}}{\inf}%
\Lambda\left(  u,\omega\right)  <\hat{c}.
\]

\end{lemma}

\textbf{Proof.} By definition of $\hat{c}$ there exists a sequence $\left(
u_{n},\omega_{n}\right)  $ in $X\times\mathbb{R}^{+}$ with $\omega_{n}\geq m$
and such that%
\[
\Lambda\left(  u_{n},\omega_{n}\right)  \rightarrow\hat{c}.
\]
Clearly $\omega_{n}$ is bounded and consequently also $\frac{J\left(
u_{n}\right)  }{K\left(  u_{n}\right)  }$ is bounded. So, up to a subsequence,
we have%
\[
\omega_{n}\rightarrow\bar{\omega}\geq m\text{ and }a_{n}\rightarrow\bar
{a},\text{ }a_{n}=\frac{J\left(  u_{n}\right)  }{K\left(  u_{n}\right)  }.
\]
Then
\[
\hat{c}=\frac{1}{2}\left(  \frac{\bar{a}}{\bar{\omega}}+\bar{\omega}\right)  .
\]
We claim that
\begin{equation}
\bar{a}<m^{2}.\label{claim}%
\end{equation}
In fact
\begin{align}
\hat{c}  & =\frac{1}{2}\left(  \bar{a}\frac{1}{\bar{\omega}}+\bar{\omega
}\right) \label{mi}\\
& =\frac{1}{2}\left(  \frac{m^{2}}{\bar{\omega}}+\bar{\omega}\right)
-\frac{1}{2}(m^{2}-\bar{a})\frac{1}{\bar{\omega}}.\nonumber
\end{align}

Then, by Lemma \ref{stupido} and (\ref{mi}), we get%
\begin{equation}
m>\frac{1}{2}\left(  \frac{m^{2}}{\bar{\omega}}+\bar{\omega}\right)  -\frac
{1}{2}(m^{2}-\bar{a})\frac{1}{\bar{\omega}}.\label{since}%
\end{equation}

On the other hand
\begin{equation}
\frac{1}{2}\left(  \frac{m^{2}}{\bar{\omega}}+\bar{\omega}\right)  \geq
m,\label{facile}%
\end{equation}

then (\ref{since}) and (\ref{facile}) imply that%
\[
-\frac{1}{2}(m^{2}-\bar{a})\frac{1}{\bar{\omega}}<0.
\]
So (\ref{claim}) is proved.

Now by (\ref{claim}) we can take $\hat{\omega}$ such that
\[
m>\hat{\omega}>\sqrt{\bar{a}},
\]
and, since $\bar{\omega}\geq m,$ we have
\[
\bar{\omega}>\hat{\omega}>\sqrt{\bar{a}}.
\]
So it can be easily deduced that%
\[
\frac{1}{2}\left(  \frac{\bar{a}}{\hat{\omega}}+\hat{\omega}\right)  <\frac
{1}{2}\left(  \frac{\bar{a}}{\bar{\omega}}+\bar{\omega}\right)  .
\]
Then%
\[
\lim\Lambda\left(  u_{n},\hat{\omega}\right)  =\frac{1}{2}\left(  \frac
{\bar{a}}{\hat{\omega}}+\hat{\omega}\right)  <\frac{1}{2}\left(  \frac{\bar
{a}}{\bar{\omega}}+\bar{\omega}\right)  =\hat{c}.
\]
So, for $n$ large, we have $\Lambda\left(  u_{n},\hat{\omega}\right)  <\hat{c}
$ and the conclusion follows.

$\square$

Now we are ready to prove Theorem \ref{abstract}.

\textbf{Proof of Th. \ref{abstract}.} By Lemma \ref{pomponia} the set $\Sigma$
defined in (\ref{defined}) is not empty. Let $\sigma\in\Sigma$ and $\left(
u_{n},\omega_{n}\right)  $ be a minimizing sequence for $E$ on $M_{\sigma}$.
By standard variational arguments (see e.g. \cite{ar}, \cite{Struwe}) we can
assume that $\left(  u_{n},\omega_{n}\right)  $ is also a P.S. sequence. Since
$\sigma\in\Sigma,$ we have%
\[
c=\lim\Lambda\left(  u_{n},\omega_{n}\right)  =\inf\left\{  \Lambda\left(
u,\omega\right)  :\left(  u,\omega\right)  \in M_{\sigma}\right\}  <\hat{c}.
\]
Then, by the lemma \ref{cinzia}, $\left(  u_{n},\omega_{n}\right)  $ possess a
strongly convergent subsequence and hence $E$ has a minimizer on $M_{\sigma}$.
Let us finally show that $\Sigma$ is open. Take $\sigma\in\Sigma;$ we have to
prove that, for $\varepsilon$ small, $\sigma+\varepsilon\in\Sigma.$ Let
$\left(  u_{0},\omega_{0}\right)  $ be a minimizer of $E$ on $M_{\sigma},$
then, since $\sigma\in\Sigma,$ we have
\begin{equation}
\Lambda\left(  u_{0},\omega_{0}\right)  <\hat{c}.\label{let}%
\end{equation}

Since $2\omega_{0}K(u_{0})=\sigma,$ by definition of $M_{\sigma+\varepsilon},$
we have%
\begin{equation}
\left(  u_{0},\omega_{0}+\frac{\varepsilon}{2K(u_{0})}\right)  \in
M_{\sigma+\varepsilon}.\label{poi}%
\end{equation}

Then%
\begin{equation}
\underset{\left(  u,\omega\right)  \in M_{\sigma+\varepsilon}}{\inf}%
\Lambda\left(  u,\omega\right)  \leq\Lambda\left(  u_{0},\omega_{0}%
+\frac{\varepsilon}{2K(u_{0})}\right)  .\label{quindi}%
\end{equation}

By (\ref{let}) and by (\ref{quindi}) we easily deduce that for $\varepsilon$
small we have%
\[
\underset{\left(  u,\omega\right)  \in M_{\sigma+\varepsilon}}{\inf}%
\Lambda\left(  u,\omega\right)  <\hat{c}.
\]

$\square$

\bigskip

\section{$Q$-balls}

\subsection{The Nonlinear Klein-Gordon equation}

In this section we will apply the abstract Theorem \ref{abstract} to the
existence of hylomorphic solitons of the nonlinear Klein-Gordon equation
(NKG):
\begin{equation}
\square\psi+W^{\prime}(\psi)=0\tag{NKG}\label{KG}%
\end{equation}
where $\square=\partial_{t}^{2}-\nabla^{2}$,$\;\psi:\mathbb{R}^{N}%
\rightarrow\mathbb{C}$ ($N\geq3$) and $W:\mathbb{C}\rightarrow\mathbb{R}$
with
\begin{equation}
W(\psi)=F(|\psi|)\label{www}%
\end{equation}
for some smooth function $F:\mathbb{R}^{+}\rightarrow\mathbb{R}$ and
\[
W^{\prime}(\psi)=F^{\prime}(\left\vert \psi\right\vert )\frac{\psi}{\left\vert
\psi\right\vert }.
\]
In particular we are interested in the existence of $Q-\mathit{balls.}$
Coleman called $Q-\mathit{balls}$ (\cite{Coleman86}) those solitary waves of
(NKG) which are spherically symmetric and this is the name generally used in
Physics literature. From now on, we always will assume that
\begin{equation}
W(0)=W^{\prime}(0)=0.\label{W0}%
\end{equation}
Eq. (\ref{KG}) is the Euler-Lagrange equation of the action functional
\begin{equation}
\int\left(  \frac{1}{2}\left\vert \partial_{t}\psi\right\vert ^{2}-\frac{1}%
{2}|\nabla\psi|^{2}-W(\psi)\right)  dxdt.\label{az}%
\end{equation}

Sometimes it will be useful to write $\psi$ in polar form, namely
\begin{equation}
\psi(t,x)=u(t,x)e^{iS(t,x)}\label{polar}%
\end{equation}
where $u(t,x)\in\mathbb{R}^{+}$ and $S(t,x)\in\mathbb{R}/(2\pi\mathbb{Z}); $
if we set $u_{t}=\partial_{t}u,$%
\begin{equation}
\mathbf{k}(t,x)=\nabla S(t,x)\label{bimbe}%
\end{equation}
and
\begin{equation}
\omega(t,x)=-\partial_{t}S(t,x),\label{belle}%
\end{equation}
the state ${\Psi}$ is uniquely defined by the quadruple $(u,u_{t}%
,\omega,\mathbf{k})$. Using these variables, the action $\mathcal{S}%
=\int\mathcal{L}dxdt$ takes the form
\begin{equation}
\mathcal{S}(u,u_{t},\omega,\mathbf{k})=\frac{1}{2}\int\left[  u_{t}%
^{2}-\left\vert \nabla u\right\vert ^{2}+\left(  \omega^{2}-\mathbf{k}%
^{2}\right)  u^{2}\right]  dxdt-\int W(u)dxdt=0
\end{equation}
and equation (\ref{KG}) becomes:%

\begin{equation}
\square u+\left(  \mathbf{k}^{2}-\omega^{2}\right)  u+W^{\prime}%
(u)=0\label{KG1}%
\end{equation}

\begin{equation}
\partial_{t}\left(  \omega u^{2}\right)  +\nabla\cdot\left(  \mathbf{k}%
u^{2}\right)  =0.\label{KG2}%
\end{equation}

\ 

The energy and the charge take the following form:
\begin{equation}
E(\Psi)=\int\left[  \frac{1}{2}\left\vert \partial_{t}\psi\right\vert
^{2}+\frac{1}{2}\left\vert \nabla\psi\right\vert ^{2}+W(\psi)\right]
dx\label{energy}%
\end{equation}%
\begin{equation}
H(\Psi)=-\operatorname{Im}\int\partial_{t}\psi\overline{\psi}\;dx.\label{im}%
\end{equation}

(the sign "minus"in front of the integral is a useful convention).

Using (\ref{polar}) we get:
\begin{equation}
E(u,u_{t},\omega,\mathbf{k})=\int\left[  \frac{1}{2}\left(  \partial
_{t}u\right)  ^{2}+\frac{1}{2}\left\vert \nabla u\right\vert ^{2}+\frac{1}%
{2}\left[  \omega^{2}+\mathbf{k}^{2}\right]  u^{2}+W(u)\right]
dx\label{penergy}%
\end{equation}%
\begin{equation}
H(u,\omega)=\int\omega\,u^{2}dx.\label{cha}%
\end{equation}

A particular type of solutions of eq. (\ref{KG}) are the \textit{standing
waves}. A \textit{standing wave} is a finite energy solution of (\ref{KG})
having the following form
\begin{equation}
\psi_{0}(t,x)=u(x)e^{-i\omega t}\text{, }u\text{ }\geq0.\label{sw}%
\end{equation}

\noindent Substituting (\ref{sw}) in eq. (\ref{KG}), we get
\begin{equation}
-\Delta u+W^{\prime}(u)=\omega^{2}u,\text{ }u\geq0.\label{static}%
\end{equation}

Let $N=3.$ Since the action functional (\ref{az}) is invariant for the Lorentz
group, we can obtain other solutions $\psi_{\mathbf{v}}(t,x)$ just making a
Lorentz transformation on it. Namely, if we take the velocity $\mathbf{v}%
=(v,0,0),$ $\left\vert v\right\vert <1$, and set
\[
t^{\prime}=\gamma\left(  t-vx_{1}\right)  ,\text{ }x_{1}^{\prime}%
=\gamma\left(  x_{1}-vt\right)  ,\text{ }x_{2}^{\prime}=x_{2},\text{ }%
x_{3}^{\prime}=x_{3}\;\;\;\text{with}\;\;\;\gamma=\frac{1}{\sqrt{1-v^{2}}},
\]
it turns out that
\[
\psi_{\mathbf{v}}(t,x)=\psi(t^{\prime},x^{\prime})
\]
is a solution of (\ref{KG}).

More exactly, given a standing wave $\psi(t,x)=u(x)e^{-i\omega t},$ the
function $\psi_{\mathbf{v}}(t,x):=\psi(t^{\prime},x^{\prime})$ is a solitary
wave which travels with velocity $\mathbf{v.}$ Thus, if $u(x)=u(x_{1}%
,x_{2},x_{3})$ is any solution of Eq. (\ref{static}), then
\begin{equation}
\psi_{\mathbf{v}}(t,x_{1},x_{2},x_{3})=u\left(  \gamma\left(  x_{1}-vt\right)
,x_{2},x_{3}\right)  e^{i(\mathbf{k}_{\mathbf{v}}\mathbf{\cdot x}%
-\omega_{\mathbf{v}}t)}\;\text{ }\label{solitone}%
\end{equation}
is a solution of Eq. (\ref{KG}) provided that%

\begin{equation}
\omega_{\mathbf{v}}=\gamma\omega\;\;\text{and\ \ }\;\mathbf{k}_{\mathbf{v}%
}=\gamma\omega\mathbf{v.}\label{kome}%
\end{equation}

\subsection{Existence results for Q-balls\label{qball}}

We write $W$ as follows%
\begin{equation}
W(s)=\frac{m^{2}}{2}s^{2}+N(s),\text{ }s\geq0;\label{NN}%
\end{equation}
and we will identify $W(s)$ with $F(s).$ We make the following assumptions:

\begin{itemize}
\item (W-i) \textbf{(Positivity}) $W(s)\geq0$

\item (W-ii) \textbf{(Nondegeneracy}) $W=$ $W(s)$ ( $s\geq0)$ is $C^{2}$ near
the origin with $W(0)=W^{\prime}(0)=0;\;W^{\prime\prime}(0)=m^{2}$\ $>0$

\item (W-iii) \textbf{(Hylomorphy}) $\exists s_{0}:\;N(s_{0})<0$\ 

\item (W-iiii)\textbf{(Growth condition}) Al least one of the following
assumptions holds:

\begin{itemize}
\item (a) there are constants $a,b>0,$ $2<p<2N/(N-2)$ such that for any $s>0:
$%
\[
|N^{\prime}(s)|\ \leq as^{p-1}+bs^{2-\frac{2}{p}}.
\]

\item (b) $\exists s_{1}>s_{0}:$ $N^{\prime}(s_{1})\geq0.$
\end{itemize}
\end{itemize}

Here there are some comments on assumptions (W-i), (W-ii), (W-iii), (W-iiii).

(W-i) As we shall see (W-i) implies that the energy is positive; if this
condition does not hold, it is possible to have solitary waves, but not
hylomorphic waves (cf. Proposition.\ref{nonex}).

(W-ii) In order to have solitary waves it is necessary to have $W^{\prime
\prime}(0)\geq0.$ There are some results also when $W^{\prime\prime}(0)=0 $
(null-mass case, see e.g. \cite{Beres-Lions} and \cite{BBR07}), however the
most interesting situation occurs when $W^{\prime\prime}(0)>0.$

(W-iii) This is the crucial assumption which characterizes the potentials
which might produce hylomorphic solitons. As we will see, this assumption
permits to have states $\Psi$ with hylomorphy ratio $\Lambda\left(
\Psi\right)  <m$.

(W-iiii)(a) This assumption contains the usual growth condition at infinity
which guarantees the $C^{1}$ regularity of the functional. Moreover it implies
that $\left\vert N^{\prime}(s)\right\vert $ $=O($ $s^{2-\frac{2}{p}})$ for $s$ small.

If we assume alternatively (W-iiii)(b), the growth condition (W-iiii)(a) can
be avoided by using standard tricks (see Appendix).

We have the following result:

\begin{theorem}
\label{main-theorem} If (W-i),(W-ii),(W-iii),(W-iiii) hold, then there exists
an open set $\Sigma$ such that for any $\sigma\in\Sigma,$ (NKG) has a
hylomorphic soliton (see Definition\ref{hys}) of charge $\sigma$ and having
the form (\ref{sw}).
\end{theorem}

Theorem \ref{main-theorem}, in the form given here, is a very recent result
\cite{BBBM}. In fact in \cite{BBBM} it has been proved the orbital stability
of (\ref{sw}) with respect to the standard topology of $\mathcal{X}%
=H^{1}(\mathbb{R}^{N},\mathbb{C})\times L^{2}(\mathbb{R}^{N},\mathbb{C})$ and
for all the $W^{\prime}s$ which satisfy (W-i), (W-ii), (W-iii) (W-iiii).
Nevertheless Theorem \ref{main-theorem} has a very long history starting with
the pioneering paper of Rosen \cite{rosen68}. Coleman \cite{Coleman78} and
Strauss \cite{strauss} gave the first rigorous proofs of existence of
solutions of the type (\ref{sw}) for (NKG) and for some particular $W^{\prime
}s$.\ Later very general existence conditions have been found by Berestycki
and Lions \cite{Beres-Lions}. In particular, if $W$ satisfies (W-i), (W-ii),
(W-iii), (W-iiii), from their paper we can deduce (see \cite{befogranas}) the
existence of $Q$-balls of type (\ref{sw} ) for any $\omega\in\left(
\omega_{0},m\right)  $ where
\[
\omega_{0}:=\inf\left\{  \lambda>0:\,W\left(  u\right)  <\textstyle{\frac
{1}{2}}\lambda^{2}u^{2}~\text{for some }u>0\right\}  .
\]

Notice that the hylomorphy condition (W-iii) guarantees that $\omega_{0}<m,$
and hence that $\left(  \omega_{0},m\right)  \neq\varnothing.$

The first orbital stability results are due to Shatah: in \cite{shatah} a
condition for orbital stability is given; however this condition is difficult
to be verified in concrete situations. More recently \cite{BBBM} a sufficient
and (essentially) necessary condition for the orbital stability has been
proved. This condition is given directly on $W$ and it permits to deduce
immediately Theorem \ref{main-theorem}.

\ Here we study the equation (\ref{static}) with $0<\omega^{2}<m^{2}$ by using
theorem \ref{abstract} and prove a weaker version of Theorem
\ref{main-theorem}, namely we do not prove the orbital stability but we
confine ourselves to show the existence of hylomorphic waves (see Definition
\ref{hys}) for (NKG) .

In this case we set:%

\[
X=H_{r}^{1}=\left\{  u\in H^{1}(\mathbb{R}^{N}):u\ \text{is radially
symmetric}\right\}  ,
\]%
\begin{equation}
\left\langle L_{1}u,u\right\rangle =\int\left(  \left\vert \nabla u\right\vert
^{2}+m^{2}u^{2}\right)  dx;\ N_{1}(u)=\int N(u)dx,\label{t}%
\end{equation}%
\begin{align}
J(u)  & =\frac{1}{2}\left\langle L_{1}u,u\right\rangle +N_{1}(u)\label{tt}\\
& =\frac{1}{2}\int\left(  \left\vert \nabla u\right\vert ^{2}+m^{2}%
u^{2}\right)  dx+\int N(u)dx,\nonumber
\end{align}

\begin{equation}
\left\langle L_{0}u,u\right\rangle =K(u)=\frac{1}{2}\int u^{2}dx;\ N_{0}%
(u)=0.\label{ttt}%
\end{equation}

First of all we observe that by (W-iiii)(a) \ the functional $J$ is $C^{1}$.
Whereas, if assumption (W-iiii)(b) holds, our problem can be transformed in an
equivalent one for which the functional $J$ is $C^{1}$ (see Appendix). Now in
order to use Theorem \ref{abstract}, we need to prove that assumptions
(H1,2,3) and (\ref{HH}) are satisfied.

\begin{lemma}
\label{in}The functionals $J,$ $N_{i}$ $(i=0,1)$ and $K$ defined in (\ref{t}),
(\ref{tt}) and (\ref{ttt}) satisfy the assumptions (H1,2,3).
\end{lemma}

\textbf{Proof.} Clearly (H3) holds. Let us now prove that (H1) holds. Let
$u_{n}$ be a sequence in $X$ such that $J(u_{n})$ is bounded. Then, since
$W\geq0,$ we have that
\begin{equation}
\int W(u_{n})\text{ and }\int\left\vert \nabla u_{n}\right\vert ^{2}\text{ are
bounded.}\label{lim}%
\end{equation}
So in order to show that $u_{n}$ is bounded in $X$ we need to prove that%
\begin{equation}
\left\Vert u_{n}\right\Vert _{L^{2}}\text{ is bounded.}\label{bouu}%
\end{equation}
Let
\[
2^{\ast}=\frac{2N}{N-2}%
\]
denote, as usual, the critical Sobolev exponent.

By (\ref{lim}) we have that
\begin{equation}
\int\left\vert u_{n}\right\vert ^{2^{\ast}}\text{ is bounded.}\label{limm}%
\end{equation}

Let $\varepsilon>0$ and set
\[
\Omega_{n}=\left\{  x\in\mathbb{R}^{N}:\left\vert u_{n}(x)\right\vert
>\varepsilon\right\}  \text{ and }\Omega_{n}^{c}=\mathbb{R}^{N}\backslash
\Omega_{n}.
\]
By (\ref{lim}) and since $W\geq0,$ we have
\begin{equation}
\int_{\text{ }\Omega_{n}^{c}}W(u_{n})\text{ is bounded .}\label{llim}%
\end{equation}
By $W_{2})$ we can write
\[
W(s)=\frac{1}{2}s^{2}+\circ(s^{2})\text{.}%
\]
Then, if $\varepsilon$ is small enough, there is a constant $c>0$ such that
\begin{equation}
\int_{\text{ }\Omega_{n}^{c}}W(u_{n})\geq c\int_{\Omega_{n}^{c}}u_{n}%
^{2}.\label{ma}%
\end{equation}
By (\ref{llim}) and (\ref{ma}) we get that
\begin{equation}
\int_{\Omega_{n}^{c}}u_{n}^{2}\text{ is bounded.}\label{pi}%
\end{equation}
On the other hand
\begin{equation}
\int_{\Omega_{n}}u_{n}^{2}\leq\left(  \int_{\Omega_{n}}\left\vert
u_{n}\right\vert ^{2^{\ast}}\right)  ^{\frac{N-2}{N}}meas(\Omega_{n}%
)^{\frac{2}{N}}.\label{u}%
\end{equation}
By (\ref{limm}) we have that
\begin{equation}
meas(\Omega_{n})\text{ is bounded.}\label{uu}%
\end{equation}
By (\ref{u}), (\ref{uu}), (\ref{limm}) we get that
\begin{equation}
\int_{\Omega_{n}}u_{n}^{2}\text{ is bounded.}\label{ff}%
\end{equation}
So (\ref{bouu}) follows from (\ref{pi}) and (\ref{ff}).

Let us finally prove that (H2) is satisfied.

Let $\left\{  u_{n}\right\}  \subset H_{r}^{1}$
\[
u_{n}\rightharpoonup u\text{ weakly in }H_{r}^{1}.
\]

Now we distinguish two cases:

Assume first that (W-iiii)(a) holds.

Since $H_{r}^{1}$ is compactly embedded into $L^{p}(\mathbb{R}^{N})$,
$2<p<2^{\ast},$ (see \cite{Beres-Lions})$,$ we have that
\begin{equation}%
%TCIMACRO{\dint }%
%BeginExpansion
{\displaystyle\int}
%EndExpansion
\left\vert u_{n}-u\right\vert ^{p}dx\rightarrow0.\label{camp}%
\end{equation}
Now
\begin{align}
& \left\vert
%TCIMACRO{\dint }%
%BeginExpansion
{\displaystyle\int}
%EndExpansion
\left(  N^{\prime}(u_{n})-N^{\prime}(u)\right)  \left(  u_{n}-u\right)
dx\right\vert \nonumber\\
& \leq\left(
%TCIMACRO{\dint }%
%BeginExpansion
{\displaystyle\int}
%EndExpansion
\left\vert N^{\prime}(u_{n})-N^{\prime}(u)\right\vert ^{p^{\prime}}dx\right)
^{\frac{1}{p^{\prime}}}\left(
%TCIMACRO{\dint }%
%BeginExpansion
{\displaystyle\int}
%EndExpansion
\left\vert u_{n}-u\right\vert ^{p}dx\right)  ^{\frac{1}{p}},\text{ }p^{\prime
}=\frac{p}{p-1}\label{cimp}%
\end{align}
The sequence $u_{n}$ is bounded in $L^{p}(\mathbb{R}^{N})$ and in
$L^{2}(\mathbb{R}^{N}).$ So, by using (W-iiii)a, we deduce that $N^{\prime
}(u_{n})$ is bounded in $L^{p^{\prime}}(\mathbb{R}^{N}).$ Then, by
(\ref{camp}) and (\ref{cimp}), we deduce that $N^{\prime}$ satisfies
(\ref{stare}).

Finally we assume that (W-iiii)(b) holds.

Clearly%
\begin{equation}
u_{n}\rightarrow u\text{ strongly in }L^{p}(B_{R})\label{a}%
\end{equation}

where $R>0$ and
\[
B_{R}=\left\{  x\in\mathbb{R}^{N}:\left\vert x\right\vert <R\right\}  .
\]

Since we can assume $N^{\prime}(s)$ linear for large $s$ (see Appendix), we have%

\begin{equation}
N^{\prime}(u_{n})\rightarrow N^{\prime}(u)\text{ in }L^{2}(B_{R}).\label{c}%
\end{equation}
Now%
\begin{equation}%
%TCIMACRO{\dint }%
%BeginExpansion
{\displaystyle\int}
%EndExpansion
\left\vert N^{\prime}(u_{n})-N^{\prime}(u)\right\vert ^{2}dx=%
%TCIMACRO{\dint _{B_{R}}}%
%BeginExpansion
{\displaystyle\int_{B_{R}}}
%EndExpansion
\left\vert N^{\prime}(u_{n})-N^{\prime}(u)\right\vert ^{2}dx+%
%TCIMACRO{\dint _{B_{R}^{c}}}%
%BeginExpansion
{\displaystyle\int_{B_{R}^{c}}}
%EndExpansion
\left\vert N^{\prime}(u_{n})-N^{\prime}(u)\right\vert ^{2}dx\label{da}%
\end{equation}

and%
\begin{equation}%
%TCIMACRO{\dint _{B_{R}^{c}}}%
%BeginExpansion
{\displaystyle\int_{B_{R}^{c}}}
%EndExpansion
\left\vert N^{\prime}(u_{n})-N^{\prime}(u)\right\vert ^{2}dx=%
%TCIMACRO{\dint _{B_{R}^{c}}}%
%BeginExpansion
{\displaystyle\int_{B_{R}^{c}}}
%EndExpansion
\left\vert N^{^{\prime\prime}}(\xi_{n})\right\vert ^{2}\left\vert
u_{n}-u\right\vert ^{2}dx\label{daa}%
\end{equation}
where%
\[
B_{R}^{c}=\mathbb{R}^{N}-B_{R}%
\]%
\[
\xi_{n}(x)=tu_{n}(x)+(1-t)u(x),\text{ }0\leq t\leq1.
\]
In the following $c_{1},c_{2},c_{3}$ will denote positive constants. By a well
known radial lemma \cite{Beres-Lions} and since $\left\Vert u_{n}\right\Vert
_{X}$ is bounded, we have that for $\left\vert x\right\vert $ large
\begin{equation}
\left\vert \xi_{n}(x)\right\vert \leq\left\vert u(x)\right\vert +\left\vert
u_{n}(x)\right\vert \leq c_{1}\frac{\left\Vert u\right\Vert _{X}+\left\Vert
u_{n}\right\Vert _{X}}{\left\vert x\right\vert ^{\frac{N-1}{2}}}\leq
\frac{c_{2}}{\left\vert x\right\vert ^{\frac{N-1}{2}}}.\label{rad}%
\end{equation}

Let $\varepsilon>0$, since $N^{\prime\prime}$ is continuous in $0$ and
$N^{\prime\prime}(0)=0,$ we have, by using (\ref{rad}), that
\begin{equation}
\left\vert N^{^{\prime\prime}}(\xi_{n}(x))\right\vert ^{2}<\varepsilon\text{
for }\left\vert x\right\vert >R,\text{ }R\text{ large.}\label{do}%
\end{equation}
So, by (\ref{daa}) and (\ref{do}) and since $\left\Vert u_{n}\right\Vert
_{L^{2}}$ is bounded, we get%
\begin{equation}%
%TCIMACRO{\dint _{B_{R}^{c}}}%
%BeginExpansion
{\displaystyle\int_{B_{R}^{c}}}
%EndExpansion
\left\vert N^{\prime}(u_{n})-N^{\prime}(u)\right\vert ^{2}dx<\varepsilon%
%TCIMACRO{\dint _{B_{R}^{c}}}%
%BeginExpansion
{\displaystyle\int_{B_{R}^{c}}}
%EndExpansion
\left\vert u_{n}-u\right\vert ^{2}dx\leq\varepsilon c_{3}.\label{di}%
\end{equation}
Then by (\ref{da}), (\ref{di}) we have
\begin{equation}%
%TCIMACRO{\dint }%
%BeginExpansion
{\displaystyle\int}
%EndExpansion
\left\vert N^{\prime}(u_{n})-N^{\prime}(u)\right\vert ^{2}dx\leq\varepsilon
c_{3}+%
%TCIMACRO{\dint _{B_{R}}}%
%BeginExpansion
{\displaystyle\int_{B_{R}}}
%EndExpansion
\left\vert N^{\prime}(u_{n})-N^{\prime}(u)\right\vert ^{2}dx.\label{de}%
\end{equation}

So by (\ref{c}) and (\ref{de}) we get%

\[
N^{\prime}(u_{n})\rightarrow N^{\prime}(u)\text{ strongly in }L^{2}%
(\mathbb{R}^{N}).
\]
Then $N$ satisfies (\ref{stare}).

$\square$

\begin{lemma}
\label{seguente}Assumption (\ref{HH}) is satisfied.
\end{lemma}

\textbf{Proof.} Let $R>0$ and consider the map $u_{R}$ defined as follows
\begin{equation}
u_{R}(x)=\left\{
\begin{array}
[c]{cc}%
s_{0} & if\;\;|x|<R\\
0 & if\;\;|x|>R+1\\
s_{0}\left(  1+R-|x|\right)  & if\;\;R\leq|x|\leq R+1
\end{array}
\right. \label{mapp}%
\end{equation}

where $s_{0}$ is a such that $N(s_{0})<0.$

Clearly%
\[
\frac{J(u_{R})}{K(u_{R})}=\frac{%
%TCIMACRO{\dint }%
%BeginExpansion
{\displaystyle\int}
%EndExpansion
\left\vert \nabla u_{R}\right\vert ^{2}}{\frac{1}{2}%
%TCIMACRO{\dint }%
%BeginExpansion
{\displaystyle\int}
%EndExpansion
u_{R}^{2}}+m^{2}+\frac{%
%TCIMACRO{\dint }%
%BeginExpansion
{\displaystyle\int}
%EndExpansion
N(u_{R})}{\frac{1}{2}%
%TCIMACRO{\dint }%
%BeginExpansion
{\displaystyle\int}
%EndExpansion
u_{R}^{2}}.
\]

Easy estimates show that for $R$ large
\begin{align}%
%TCIMACRO{\dint }%
%BeginExpansion
{\displaystyle\int}
%EndExpansion
\left\vert \nabla u_{R}\right\vert ^{2}  & \leq c_{0}R^{N-1}\label{bag}\\
c_{2}R^{N}  & \leq\frac{1}{2}\int u_{R}^{2}dx\leq c_{1}R^{N}\label{beg}\\
\int N(u_{R})dr  & \leq N(s_{0})R^{N}+c_{3}R^{N-1}\label{big}%
\end{align}

where $c_{0},...,c_{3}$ are positive constants.

Then for $R$ large, since $N(s_{0})<0,$ we have%
\[
\frac{J(u_{R})}{K(u_{R})}\leq\frac{c_{0}}{c_{2}}\frac{1}{R}+m^{2}%
+\frac{N(s_{0})R^{N}}{c_{1}R^{N}}+\frac{c_{3}R^{N-1}}{c_{2}R^{N}}<m^{2}.
\]

$\square$

Assumption (W-i) is a necessary condition for the existence of hylomorphic
waves (Definition \ref{hys}), in fact the following proposition holds:.

\begin{proposition}
\label{nonex}If (W-i) does not hold, then for any $\sigma>0,$ $E(u)$ is not
bounded from below on $M_{\sigma}$.
\end{proposition}

\textbf{Proof. }Let $\sigma>0$ and assume that there exists $s_{0}$ is a such
that $W(s_{0})<0.$ We set $\Psi_{R}=(u_{R},-i\omega_{R}u_{R})$\ where $u_{R}$
is defined in (\ref{mapp}) and%
\[
\omega_{R}=\frac{\sigma}{\int u_{R}^{2}dx}.
\]

Clearly%
\begin{equation}
\omega_{R}=\frac{\sigma}{\int u_{R}^{2}dx}\leq c_{4}R^{-N}.\label{bog}%
\end{equation}

Then by (\ref{bag}), (\ref{beg}), (\ref{big}) (where $W$ replaces $N$) we have%
\begin{align*}
E\left(  \Psi_{R}\right)   & =\int\left[  \frac{1}{2}\left\vert \nabla
u_{R}\right\vert ^{2}+W(u_{R})\right]  dx+\frac{1}{2}\omega_{R}^{2}\int
u_{R}^{2}dx\\
& =\int\left[  \frac{1}{2}\left\vert \nabla u_{R}\right\vert ^{2}%
+W(u_{R})\right]  dx+\frac{1}{2}\omega_{R}\sigma\\
& \leq\frac{1}{2}c_{0}R^{N-1}+W(s_{0})R^{N}+c_{3}R^{N-1}+c_{5}R^{-N}.
\end{align*}
Hence%
\[
\underset{R\rightarrow\infty}{\lim}E\left(  \Psi_{R}\right)  =-\infty
\]

$\square$

\begin{remark}
\label{condizioni-necessarie} If (W-i) is violated, it is still possible to
have orbitally stable solitary waves (see \cite{shatah}) which are only local
minimizers. They can be destroyed by a perturbation which send them out of the
basin of attraction and are not considered solitons according to Def.
\ref{hys}.
\end{remark}

\begin{remark}
\label{puffo}We observe that the constant $m$ defined by (W-ii) coincides with
the constant $m$ defined by (\ref{brutta}) and the constant $m$ defined by
(\ref{norma}).
\end{remark}

\section{Vortices}

\subsection{Main features}

A ($hylomorphic)$ $vortex$ is a (hylomorphic) solitary wave with nonvanishing
\textit{angular momentum}. The angular momentum, by definition, is the
quantity which is preserved by virtue of the invariance under space rotations
(with respect to the origin) of the Lagrangian (see e.g.\cite{Gelfand}). In
this section we shall analyze elementary properties of the angular momentum
for (\ref{KG}) in three space dimensions; of course, making obvious changes,
the analysis includes also the two dimensional case .

The angular momentum for the solutions of (\ref{KG}) is given by
\begin{equation}
\mathbf{M}(\Psi)=\operatorname{Re}\int\mathbf{x}\times\nabla\psi\left(
\overline{\partial_{t}\psi}\right)  \;dx.
\end{equation}
Using the polar form (\ref{polar}), it can be written
\begin{equation}
\mathbf{M}(\Psi)=\int\left(  \mathbf{x}\times\nabla S\left(  \partial
_{t}Su^{2}\right)  +\mathbf{x}\times\nabla u\,\left(  \partial_{t}u\right)
\right)  \;dx.\label{amom}%
\end{equation}
where $\times$ denotes the wedge product.

It is immediate to check that standing waves (\ref{sw}) have $\mathbf{M}%
\left(  \Psi\right)  =\mathbf{0.}$ However, if we consider:
\begin{equation}
\psi\left(  t,x\right)  =\psi_{0}\left(  x\right)  e^{-i\omega t}%
\,,\quad\omega>0\,\label{stationary}%
\end{equation}
where $\psi_{0}\left(  x\right)  $ is allowed to have complex values, it is
possible to have $\mathbf{M}\left(  \Psi\right)  \neq\mathbf{0.}$ Thus, we are
led to make an ansaz of the following form:
\begin{equation}
\psi\left(  t,x\right)  =u\left(  x\right)  e^{i\left(  \ell\theta\left(
x\right)  -\omega t\right)  }\,,\quad u\left(  x\right)  \geq0,~\omega
\in\mathbb{R},\;\ell\in\mathbb{Z}-\left\{  0\right\} \label{ansatz}%
\end{equation}
and
\[
\theta\left(  x\right)  =\operatorname{Im}\log(x_{1}+ix_{2})\in\mathbb{R}%
/2\pi\mathbb{Z};\mathbb{\;\;}x=(x_{1},x_{2},x_{3}).
\]
Moreover, we assume that%
\begin{equation}
u(x)=u(r,x_{3}),\ \text{where }r=\sqrt{x_{1}^{2}+x_{2}^{2}}.\label{mariella}%
\end{equation}

By this ansaz, equation (NKG) (in the form (\ref{KG1}), (\ref{KG2})) is
equivalent to the system
\[
\left\{
\begin{array}
[c]{l}%
-\triangle u+\ell^{2}\left\vert \nabla\theta\right\vert ^{2}u+W^{\prime
}\left(  u\right)  =\omega^{2}u\\
u\triangle\theta+2\nabla u\cdot\nabla\theta=0\,.
\end{array}
\right.
\]
By the definition of $\theta$ and (\ref{mariella}) we have
\[
\triangle\theta=0\,,\quad\nabla\theta\cdot\nabla u=0\,,\quad\left\vert
\nabla\theta\right\vert ^{2}=\frac{1}{r^{2}}.
\]
where the dot $\cdot$ denotes the euclidean scalar product.

So the above system reduces to
\begin{equation}
-\triangle u+\frac{\ell^{2}}{r^{2}}u+W^{\prime}\left(  u\right)  =\omega
^{2}u\qquad\text{in }\mathbb{R}^{3}.\label{elliptic eq}%
\end{equation}
Direct computations show that the energy (\ref{energy}), the angular momentum
(\ref{amom}) and the hylenic charge (\ref{im}) become
\begin{equation}
E\left(  u\left(  x\right)  e^{i\left(  \ell\theta\left(  x\right)  -\omega
t\right)  }\right)  =\int_{\mathbb{R}^{3}}\left[  \frac{1}{2}\left\vert \nabla
u\right\vert ^{2}+\frac{1}{2}\left(  \frac{\ell^{2}}{r^{2}}+\omega^{2}\right)
u^{2}+W\left(  u\right)  \right]  dx\label{energy finite}%
\end{equation}%
\begin{equation}
\mathbf{M}\left(  u\left(  x\right)  e^{i\left(  \ell\theta\left(  x\right)
-\omega t\right)  }\right)  =-\left(  0,0,\omega\ell\int_{\mathbb{R}^{3}}%
u^{2}dx\right)  .\label{momentum nonzero}%
\end{equation}%
\begin{equation}
H\left(  u\left(  x\right)  e^{i\left(  \ell\theta\left(  x\right)  -\omega
t\right)  }\right)  =\int\omega\,u^{2}dx.
\end{equation}

The existence of vortices is an interesting and old issue in many questions of
mathematical physics as superconductivity, classical and quantum field theory,
string and elementary particle theory (see the pioneering papers \cite{ab},
\cite{nil} and e.g. the more recent ones \cite{Kim93}, \cite{Volk},
\cite{Volk-Wohn}, \cite{vil}, \cite{spinning review} with their references).

From mathematical viewpoint, the existence of vortices for (NKG) and for
(NKGM) has been studied in some recent papers ( \cite{Be-Visc}, \cite{BBR07},
\cite{bbr08}, \cite{cossu}, \cite{befov07}, \cite{befo08}).

\subsection{Existence of two dimensional vortices}

In this paper we want to apply theorem \ref{abstract} to the study of
vortices; this is possible for $N=2.$ We get the following theorem:

\begin{theorem}
\label{THM: W} Let $W:\mathbb{C}\rightarrow\mathbb{R}$ satisfy (W-i), (W-ii),
(W-iii), (W-iiii) and fix $\ell\in\mathbb{Z}-\left\{  0\right\}  $;\ then
there exists an open set $\Sigma$ such that for any $\sigma\in\Sigma,$
equation NKG has a hylomorphic vortex of the form (\ref{ansatz}).
\end{theorem}

In this case we set:%

\[
\left\langle L_{1}u,u\right\rangle =\int\left[  \left\vert \nabla u\right\vert
^{2}+\left(  \frac{\ell^{2}}{r^{2}}+m^{2}\right)  u^{2}\right]  dx;\ N_{1}%
(u)=\int N(u)dx
\]%
\[
X=\left\{  u\in H^{1}(\mathbb{R}^{2}):u\ \text{is radially symmetric\ and\ }%
\left\langle L_{1}u,u\right\rangle <\infty\right\}
\]%
\begin{align*}
J(u)  & =\frac{1}{2}\left\langle L_{1}u,u\right\rangle +N_{1}(u)\\
& =\frac{1}{2}\int\left[  \left\vert \nabla u\right\vert ^{2}+\left(
\frac{\ell^{2}}{r^{2}}+m^{2}\right)  u^{2}\right]  dx+\int N(u)dx
\end{align*}%
\[
\left\langle L_{0}u,u\right\rangle =K(u)=\frac{1}{2}\int u^{2}dx;\ N_{0}(u)=0.
\]

\begin{lemma}
\label{v1} Assumptions (H1), (H2), (H3) are satisfied
\end{lemma}

Proof. Clearly assumption (H3) is satisfied. Let us prove that assumption (H1)
is satisfied.

Let $u_{n}$ be a sequence in $X$ such that $J(u_{n})$ is bounded. Then clearly
also the sequences%
\begin{equation}%
%TCIMACRO{\dint }%
%BeginExpansion
{\displaystyle\int}
%EndExpansion
\left\vert \nabla u_{n}\right\vert ^{2},%
%TCIMACRO{\dint }%
%BeginExpansion
{\displaystyle\int}
%EndExpansion
\frac{u_{n}^{2}}{r^{2}},%
%TCIMACRO{\dint }%
%BeginExpansion
{\displaystyle\int}
%EndExpansion
W(u_{n})\label{31}%
\end{equation}

are bounded. We have to show that $u_{n}$ is bounded in $L^{2}.$ Let us first
show that there exists $M_{1}$ such that for all $n$
\begin{equation}
\left\Vert u_{n}\right\Vert _{L^{\infty}}\leq M_{1}.\label{infinity}%
\end{equation}
In fact for $u\in C_{0}^{\infty}(\mathbb{R}^{2}\backslash0),$ $u$ radially
symmetric, we set $u(x)=v(r)$ $r=\left\vert x\right\vert ,$ then%
\begin{align}
\frac{1}{2}u^{2}(x)  & =\frac{1}{2}v(r)^{2}=%
%TCIMACRO{\dint _{+\infty}^{r}}%
%BeginExpansion
{\displaystyle\int_{+\infty}^{r}}
%EndExpansion
v(r)v^{\prime}(r)dr\leq\nonumber\\
\left(
%TCIMACRO{\dint _{0}^{+\infty}}%
%BeginExpansion
{\displaystyle\int_{0}^{+\infty}}
%EndExpansion
\frac{v(r)^{2}}{r}dr%
%TCIMACRO{\dint _{0}^{+\infty}}%
%BeginExpansion
{\displaystyle\int_{0}^{+\infty}}
%EndExpansion
v^{\prime}(r)^{2}rdr\right)  ^{\frac{1}{2}}  & \leq c_{1}\left(
%TCIMACRO{\dint _{\mathbb{R}^{2}}}%
%BeginExpansion
{\displaystyle\int_{\mathbb{R}^{2}}}
%EndExpansion
\frac{u^{2}}{r^{2}}dx%
%TCIMACRO{\dint _{\mathbb{R}^{2}}}%
%BeginExpansion
{\displaystyle\int_{\mathbb{R}^{2}}}
%EndExpansion
\left\vert \nabla u\right\vert ^{2}dx\right)  ^{\frac{1}{2}}\label{32}%
\end{align}

Then, since the sequences (\ref{31}) are bounded, by (\ref{32}) we get
(\ref{infinity}).

Let $\varepsilon>0$ and set
\[
\Omega_{n}=\left\{  x\in\mathbb{R}^{2}:\left\vert u_{n}(x)\right\vert
>\varepsilon\right\}  \text{ and }\Omega_{n}^{c}=\mathbb{R}^{2}\backslash
\Omega_{n}.
\]
Then, by (\ref{infinity}), we have%
\begin{align}
\int_{\text{ }\Omega_{n}}u_{n}^{2}  & \leq\left(  \int_{\text{ }\Omega_{n}%
}u_{n}^{6}\right)  ^{\frac{1}{3}}\left(  meas(\Omega_{n})\right)  ^{\frac
{2}{3}}\leq\label{conf}\\
& \leq\left\Vert u_{n}\right\Vert _{L^{\infty}}^{2}meas(\Omega_{n})\leq
M_{1}^{2}meas(\Omega_{n}).\nonumber
\end{align}
On the other hand, if $\varepsilon$ is small enough we have (see (\ref{ma}) in
the proof of Lemma.\ref{in})%
\begin{equation}
\int_{\text{ }\Omega_{n}^{c}}W(u_{n})\geq c_{2}\int_{\text{ }\Omega_{n}^{c}%
}u_{n}^{2}.\label{oc}%
\end{equation}

Since
\[
\int_{\text{ }}W(u_{n})\leq M_{2},
\]
by (\ref{conf}) and (\ref{oc}) we deduce%
\begin{equation}
\int_{\text{ }}u_{n}^{2}=\int_{\Omega_{n}\text{ }}u_{n}^{2}+\int_{\text{
}\Omega_{n}^{c}}u_{n}^{2}\leq M_{1}^{2}meas(\Omega_{n})+\frac{M_{2}}{c_{2}%
}.\label{otto}%
\end{equation}
Then it remains to prove that
\begin{equation}
meas(\Omega_{n})\text{ is bounded.}\label{nove}%
\end{equation}
Arguing by contradiction assume that, up to a subsequence
\begin{equation}
meas(\Omega_{n})\rightarrow\infty.\label{dieci}%
\end{equation}
By a Trudingher-Moser type inequality (see \cite{ruf} and its references) on
all $\mathbb{R}^{2},$ we have for $\alpha<4\pi$%
\begin{equation}
\int_{\text{ }}e^{\alpha u_{n}^{2}}\leq c_{3}%
%TCIMACRO{\dint }%
%BeginExpansion
{\displaystyle\int}
%EndExpansion
\left\vert \nabla u_{n}\right\vert ^{2}\text{ .}\label{ruf}%
\end{equation}

Then, taking $\alpha=1$ and since $%
%TCIMACRO{\dint }%
%BeginExpansion
{\displaystyle\int}
%EndExpansion
\left\vert \nabla u_{n}\right\vert ^{2}$ is bounded, we have%
\[
e^{\varepsilon^{2}}meas(\Omega_{n})\leq\int_{\Omega_{n}\text{ }}e^{u_{n}^{2}%
}\leq\int_{\text{ }}e^{u_{n}^{2}}\leq c_{3}%
%TCIMACRO{\dint }%
%BeginExpansion
{\displaystyle\int}
%EndExpansion
\left\vert \nabla u_{n}\right\vert ^{2}\leq M_{3}%
\]
which contradicts (\ref{dieci}).

Finally, following the same arguments used in the proof of Lemma \ref{in}, it
can be proved that also assumption (H2) is satisfied$.$

$\square$

\begin{lemma}
\label{seguente copy(1)}Assumption (\ref{HH}) is satisfied
\end{lemma}

Proof. Let $R>1$ and consider the map $u_{R}$ defined as follows
\[
u_{R}(x)=\left\{
\begin{array}
[c]{cc}%
0 & \text{if\ }\;|x|\leq R-1\text{ or }|x|\geq2R+1\\
s_{0}\left(  |x|-R+1\right)  & \text{if\ }R\geq\;|x|>R-1\\
s_{0} & 2R\geq\left\vert x\right\vert >R\\
s_{0}\left(  1+2R-|x|\right)  & \text{if\ }2R+1\geq\;|x|>2R
\end{array}
\right.
\]

where $s_{0}$ is a such that $N(s_{0})<0.$

Clearly%
\begin{equation}
\frac{J(u_{R})}{K(u_{R})}=\frac{%
%TCIMACRO{\dint }%
%BeginExpansion
{\displaystyle\int}
%EndExpansion
\left\vert \nabla u_{R}\right\vert ^{2}}{%
%TCIMACRO{\dint }%
%BeginExpansion
{\displaystyle\int}
%EndExpansion
u_{R}^{2}}+m^{2}+\frac{%
%TCIMACRO{\dint }%
%BeginExpansion
{\displaystyle\int}
%EndExpansion
\frac{\ell^{2}u_{R}^{2}}{r^{2}}}{%
%TCIMACRO{\dint }%
%BeginExpansion
{\displaystyle\int}
%EndExpansion
u_{R}^{2}}+\frac{%
%TCIMACRO{\dint }%
%BeginExpansion
{\displaystyle\int}
%EndExpansion
N(u_{R})}{\frac{1}{2}%
%TCIMACRO{\dint }%
%BeginExpansion
{\displaystyle\int}
%EndExpansion
u_{R}^{2}}.\label{clear}%
\end{equation}

Easy estimates show that for $R$ large
\begin{align*}%
%TCIMACRO{\dint }%
%BeginExpansion
{\displaystyle\int}
%EndExpansion
\left\vert \nabla u_{R}\right\vert ^{2}  & \leq c_{0}R\\%
%TCIMACRO{\dint }%
%BeginExpansion
{\displaystyle\int}
%EndExpansion
\frac{\ell^{2}u_{R}^{2}}{r^{2}}  & \leq\frac{c_{1}}{R}+c_{2}\\
\int N(u_{R})dr  & \leq c_{3}N(s_{0})R^{2}+c_{4}R\\
c_{6}R^{2}  & \geq\int u_{R}^{2}dx\geq c_{5}R^{2}%
\end{align*}
where $c_{0},...,c_{6}$ are positive constants.

Then for $R$ large, since $N(s_{0})<0,$ we have%
\[
\frac{J(u_{R})}{K(u_{R})}<m^{2}.
\]

$\square$

\section{The Nonlinear Klein-Gordon-Maxwell equations}

\subsection{General features of NKGM}

The Nonlinear Klein-Gordon-Maxwell equations (NKGM) are (see e.g.
\cite{befogranas}, \cite{befo})%

\begin{equation}
\left(  \partial_{t}+iq\varphi\right)  ^{2}\psi-\left(  \nabla-iq\mathbf{A}%
\right)  ^{2}\psi+W^{\prime}(\psi)=0\tag{NKGM-1}\label{e1}%
\end{equation}

\begin{equation}
\nabla\cdot\left(  \partial_{t}\mathbf{%
%TCIMACRO{\QTR{mathbf}{A}}%
%BeginExpansion
A%
%EndExpansion
}+\nabla\varphi\right)  =q\;\text{Im}\left(  \partial_{t}\psi\overline{\psi
}\right)  +q^{2}\varphi\left\vert \psi\right\vert ^{2}\tag{NKGM-2}\label{e2}%
\end{equation}

\begin{equation}
\nabla\times\left(  \nabla\times\mathbf{A}\right)  +\partial_{t}\left(
\partial_{t}\mathbf{%
%TCIMACRO{\QTR{mathbf}{A}}%
%BeginExpansion
A%
%EndExpansion
}+\nabla\varphi\right)  =q\;\text{Im}\left(  \nabla\psi\overline{\psi}\right)
-q^{2}\mathbf{A}\left\vert \psi\right\vert ^{2}\tag{NKGM-3}\label{e3}%
\end{equation}

where $q$ is a parameter which, in some models, is interpreted as the electron
charge and $W$ satisfies (\ref{www}). They are the Euler-Lagrange equations of
the action:
\begin{equation}
\mathcal{S}=\int\mathcal{L\;}dxdt,\;\;\mathcal{L}=\mathcal{L}_{0}%
+\mathcal{L}_{1}-W(\psi),\label{completa}%
\end{equation}
where
\begin{equation}
\mathcal{L}_{0}=\frac{1}{2}\left[  \left\vert \left(  \partial_{t}%
+iq\varphi\right)  \psi\right\vert ^{2}-\left\vert \left(  \nabla
-iq\mathbf{A}\right)  \psi\right\vert ^{2}\right]
\end{equation}%
\begin{equation}
\mathcal{L}_{1}=\frac{1}{2}\left[  \left\vert \partial_{t}\mathbf{%
%TCIMACRO{\QTR{mathbf}{A}}%
%BeginExpansion
A%
%EndExpansion
}+\nabla\varphi\right\vert ^{2}-\frac{1}{2}\left\vert \nabla\times
\mathbf{A}\right\vert ^{2}\right]  .
\end{equation}
In this case, the state of the system is given by
\[
{\Psi}=(\psi,\psi_{t},\varphi,\varphi_{t},\mathbf{A,A}_{t}).
\]

If we use the notation (\ref{polar}, \ref{bimbe}, \ref{belle}) and if we set
\begin{equation}
\mathbf{E=-}\left(  \partial_{t}\mathbf{A}+\nabla\varphi\right) \label{E}%
\end{equation}%
\begin{equation}
\mathbf{H}=\nabla\times\mathbf{A}\label{H}%
\end{equation}%
\begin{align}
\Omega & =-\left(  \partial_{t}S+q\varphi\right)  =\omega-q\varphi
\label{omega}\\
\rho & =q\Omega u^{2}%
\end{align}%
\begin{align}
\mathbf{K}  & =\nabla S-q\mathbf{A=k}-q\mathbf{A}\label{tar}\\
\mathbf{J}  & =q\mathbf{K}u^{2}.
\end{align}
Equations (\ref{e1}), (\ref{e2}), (\ref{e3}) can be written as follows (see
e.g. \cite{befogranas}):
\begin{equation}
\square u+\left(  \mathbf{K}^{2}+\Omega^{2}\right)  u+W^{\prime}%
(u)=0\tag{\textsc{matter}}\label{materia}%
\end{equation}%
\begin{equation}
\nabla\cdot\mathbf{E}=\rho\tag{\textsc{gauss}}\label{gauss}%
\end{equation}%
\begin{equation}
\nabla\times\mathbf{H}-\frac{\partial\mathbf{E}}{\partial t}=\mathbf{J}%
\tag{\textsc{ampere}}\label{ampere}%
\end{equation}
Moreover, by the positions (\ref{E}) and (\ref{H}), $\mathbf{E}$ and
$\mathbf{H}$ satisfy also the equations%
\begin{equation}
\nabla\times\mathbf{E}+\frac{\partial\mathbf{H}}{\partial t}%
=0\tag{\textsc{faraday}}\label{faraday}%
\end{equation}%
\begin{equation}
\nabla\cdot\mathbf{H}=0.\tag{\textsc{nomonopole}}\label{monopole}%
\end{equation}

The equations (\ref{gauss}),(\ref{ampere}),(\ref{faraday}),(\ref{monopole})
are the Maxwell's equations and equation (\ref{materia}) represents a model of
interaction of matter with the elecromagnetic field (see for example
\cite{befogranas}, \cite{fel} ch. 3, \cite{rub} ch. 2 in Part 1, and
\cite{yangL} ch.1).

The energy takes the following form (see \cite{befogranas}):
\[
E(\Psi)=\int\left[  \frac{1}{2}u_{t}^{2}+\frac{1}{2}\left\vert \nabla
u\right\vert ^{2}+\frac{1}{2}\left(  \mathbf{K}^{2}+\Omega^{2}\right)
u^{2}+W(u)\,+\frac{1}{2}\left(  \mathbf{E}^{2}+\mathbf{H}^{2}\right)  \right]
dx
\]
and the hylenic charge takes the form:
\[
H(\Psi)=\int\Omega u^{2}dx=\int\left(  \omega-q\varphi\right)  u^{2}.
\]
In some models, $H(\Psi),$ if positive, represents the number of particles
contained in the state $\Psi,$ otherwise, $-H(\Psi)$ represents the number of
antiparticles. The global electric charge is given by
\[
Q(\Psi)=qH(\Psi)=\int\left(  q\omega-q^{2}\varphi\right)  u^{2}.
\]
Thus, if $\psi$ is rescaled in such a way to have $q=1,$ the hylenic charge
$H(\Psi)$ and the electric charge $Q(\Psi)$ coincide.

\subsection{Existence results for the NKGM}

In this paper we are interested to apply Theorem \ref{abstract} to find
electrostatic standing waves, namely solutions of of (\ref{e1+}),(\ref{e3+}),
(\ref{e4+}), having the form
\begin{align}
\psi\left(  t,x\right)   & =u\left(  x\right)  e^{-i\omega t},\text{\ }%
u\in\mathbb{R}^{+},\ \omega\in\mathbb{R},\text{\ }s\in\frac{\mathbb{R}}%
{2\pi\mathbb{Z}}\label{static1}\\
\mathbf{A}  & =0\mathbf{,\ }\partial_{t}\varphi=0.\label{static2}%
\end{align}
Using (\ref{static1}) and (\ref{static2}), equations (\ref{e1}), (\ref{e2}),
(\ref{e3}) become:
\begin{equation}
\square u+W^{\prime}(u)+\left[  \left\vert \nabla S-q\mathbf{A}\right\vert
^{2}-\left(  \frac{\partial S}{\partial t}+q\varphi\right)  ^{2}\right]
\,u=0\label{e1+}%
\end{equation}%
\begin{equation}
\frac{\partial}{\partial t}\left[  \left(  \frac{\partial S}{\partial
t}+q\varphi\right)  u^{2}\right]  -\nabla\cdot\left[  \left(  \nabla
S-q\mathbf{A}\right)  u^{2}\right]  =0\label{e2+}%
\end{equation}%
\begin{equation}
\nabla\cdot\left(  \frac{\partial\mathbf{%
%TCIMACRO{\QTR{mathbf}{A}}%
%BeginExpansion
A%
%EndExpansion
}}{\partial t}+\nabla\varphi\right)  =q\left(  \frac{\partial S}{\partial
t}+q\varphi\right)  u^{2}\;\label{e3+}%
\end{equation}%
\begin{equation}
\nabla\times\left(  \nabla\times\mathbf{A}\right)  +\frac{\partial}{\partial
t}\left(  \frac{\partial\mathbf{%
%TCIMACRO{\QTR{mathbf}{A}}%
%BeginExpansion
A%
%EndExpansion
}}{\partial t}+\nabla\varphi\right)  =q\left(  \nabla S-q\mathbf{A}\right)
u^{2}\;.\label{e4+}%
\end{equation}

Observe that equation (\ref{e2+}) is the continuity equation
\[
\partial_{t}\rho+\nabla\cdot\mathbf{J}=0,
\]
and it easily follows from equation (\ref{e3+}) and (\ref{e4+}). Then we are
reduced to study the system (\ref{e1+}), (\ref{e3+}), (\ref{e4+}).

The existence of solitary waves for (NKGM) depends on the constant $q;$ more
exactly we have the following theorem:

\begin{theorem}
\label{main-theorem copy(1)} Assume that (W-i),(W-ii),(W-iii), (W-iiii) hold.
Then there exists a set $\Sigma_{NKGM}\subset\mathbb{R}^{2}$ such that for any
$(\sigma,q)\in\Sigma_{NKGM},$ the nonlinear Klein-Gordon -Maxwell equations
(NKGM) have an hylomorphic, electrostatic (see (\ref{static1}), (\ref{static2}%
)) wave of charge $\sigma$. Moreover $\Sigma_{NKGM}$ has the following form%
\[
\Sigma_{NKGM}=\left\{  (\sigma,q)\in\mathbb{R}^{2}:\sigma\in\Sigma
_{q},\ 0<q<q^{\ast}\right\}
\]
where $q^{\ast}>0$ and $\Sigma_{q}$ is an open set which is not empty for
$0<q<q^{\ast}.$
\end{theorem}

\begin{remark}
\bigskip The existence of electrostatic standing waves has been first analyzed
when $W(s)$ changes sign, namely when $W(s)=s^{2}-s^{p}$ $(s>0,$ $p>2)$
(\cite{bf}, \cite{ca}, \cite{tea}, \cite{tea2}). More recently also cases in
which $W\geq0$ have been considered (\cite{befo}, \cite{befo08}, \cite{mu}).
\end{remark}

If (\ref{static1}) and ( \ref{static2}) hold, equation (\ref{e4+}) is
identically satisfied, while (\ref{e1+}) and (\ref{e3+}) become
\begin{equation}
-\Delta u+W^{\prime}(u)=\left(  \omega-q\varphi\right)  ^{2}u\label{a11}%
\end{equation}%
\begin{equation}
-\Delta\varphi=q\left(  \omega-q\varphi\right)  u^{2}.\;\label{a222}%
\end{equation}

We set
\begin{equation}
\mathcal{X}_{0}=\left\{  \Psi=\left(  u(x),-i\omega u(x),\varphi
(x),0,\mathbf{0,0}\right)  ,\;u\in H^{1}(\mathbb{R}^{N}),\varphi\in
\mathcal{D}^{1,2}(\mathbb{R}^{3}),\omega\in\mathbb{R}\right\}  .
\end{equation}
Clearly $\mathcal{X}_{0}$ is a subset of the phase space which contains the
electro-static standing waves. To any state $\Psi\in\mathcal{X}_{0},$ we can
associate a triple
\[
\left(  u,\varphi,\omega\right)  \in H^{1}(\mathbb{R}^{3})\times
\mathcal{D}^{1,2}(\mathbb{R}^{3})\times\mathbb{R};
\]
the corresponding energy and charge take the following form:
\begin{equation}
E_{q}\left(  u,\varphi,\omega\right)  =\int\left[  \frac{1}{2}\left\vert
\nabla u\right\vert ^{2}+\frac{1}{2}\left\vert \nabla\varphi\right\vert
^{2}+\frac{1}{2}\Omega^{2}u^{2}+W(u)\right]  dx
\end{equation}%
\[
H_{q}\left(  u,\varphi,\omega\right)  =\int\Omega u^{2}dx
\]
where, according to (\ref{omega}),
\[
\Omega=\omega-q\varphi.
\]

Now we would like to apply theorem \ref{abstract}. Unforunately, we cannot do
it directly, since $E_{q}$ and $H_{q}$ do not satisfy the required properties,
namely they do not have the form (\ref{forma}) and (\ref{formina}). However,
we can transform this problem in such a way that Theorem \ref{abstract} can be
used. To do this, we introduce a smaller space $\mathcal{Z}_{0}\subset
\mathcal{X}_{0}$ which contains the states which satisfy equation
(\ref{a222}), namely
\begin{equation}
\mathcal{Z}_{0}=\left\{  \Psi\in\mathcal{X}_{0}:-\Delta\varphi=q\left(
\omega-q\varphi\right)  u^{2}\right\}  .\label{zed}%
\end{equation}

We remark that for $u\in H^{1}(\mathbb{R}^{3})$ and $\omega\in\mathbb{R}$
given, equation (\ref{a222}) has a unique solution $\varphi_{u}\in
\mathcal{D}^{1,2}(\mathbb{R}^{3})\;$(see \cite{bf}); then
\[
\mathcal{Z}_{0}\cong H^{1}(\mathbb{R}^{3})\times\mathbb{R}.
\]

Now we want to find a nice and useful way to write $E_{q},H_{q}$ and
$\Lambda_{q}$ restricted to $\mathcal{Z}_{0}$. First, we divide the energy in
two parts:
\begin{equation}
E_{q}\left(  u,\varphi,\omega\right)  =J\left(  u\right)  +F_{q}\left(
u,\varphi,\omega\right) \label{ener}%
\end{equation}

where
\begin{align}
J\left(  u\right)   & =\int\left[  \frac{1}{2}\left\vert \nabla u\right\vert
^{2}+W(u)\right]  dx\label{ba}\\
F_{q}\left(  u,\varphi,\omega\right)   & =\frac{1}{2}\int\left[  \left\vert
\nabla\varphi\right\vert ^{2}+\Omega^{2}u^{2}\right]  dx\label{be}%
\end{align}
Now let $u\in H^{1}(\mathbb{R}^{3})$ and consider the solution $\varphi_{u}$
of (\ref{a222}).

Multiplying both sides of equation (\ref{a222}) by $\varphi_{u}$ and
integrating, we get
\[
\int\left\vert \nabla\varphi_{u}\right\vert ^{2}dx=\int q\varphi_{u}\Omega
u^{2}.
\]
Then
\begin{align*}
F_{q}\left(  u,\varphi_{u},\omega\right)   & =\frac{1}{2}\int\left[
q\varphi_{u}\Omega u^{2}+\Omega^{2}u^{2}\right]  dx\\
& =\frac{1}{2}\omega^{2}\int\left(  1-q\frac{\varphi_{u}}{\omega}\right)
u^{2}dx.
\end{align*}
So we have
\begin{equation}
F_{q}\left(  u,\varphi_{u},\omega\right)  =\frac{1}{2}\omega^{2}\int\left(
1-q\frac{\varphi_{u}}{\omega}\right)  u^{2}dx.\label{bi}%
\end{equation}

For $u\in H^{1}(\mathbb{R}^{3}),$ let $\Phi=\Phi_{u}$ be the solution of the
equation
\begin{equation}
-\Delta\Phi_{u}+q^{2}u^{2}\Phi_{u}=qu^{2}.\label{b2}%
\end{equation}
Clearly
\begin{equation}
\varphi_{u}=\omega\Phi_{u}\label{d}%
\end{equation}
solves eq. (\ref{a222}) and we have that
\begin{equation}
F_{q}\left(  u,\varphi_{u},\omega\right)  =F_{q}\left(  u,\omega\Phi
_{u},\omega\right)  =\frac{1}{2}\,\omega^{2}\int\left(  1-q\Phi_{u}\right)
u^{2}dx=\omega^{2}K_{q}\left(  u\right)  ,\label{sea}%
\end{equation}

where
\begin{equation}
K_{q}\left(  u\right)  :=\frac{1}{2}\int\left(  1-q\Phi_{u}\right)
u^{2}dx.\label{sita}%
\end{equation}
By (\ref{ener}) and (\ref{sea}) the energy on the states contained in
$\mathcal{Z}_{0}$ (see (\ref{zed})) can be written as a functional of the two
variables $\omega$ and $u$ and having the form (\ref{forma}):%
\begin{equation}
\tilde{E}_{q}\left(  u,\omega\right)  =E_{q}\left(  u,\varphi_{u}%
,\omega\right)  =J\left(  u\right)  +\omega^{2}K_{q}\left(  u\right)
.\label{formade}%
\end{equation}
Analogously, also the hylenic charge can be expressed via the variables $u$
and $\omega$ and having the form (\ref{formina}):
\begin{align*}
\tilde{H}_{q}\left(  u,\omega\right)   & =H_{q}\left(  u,\varphi_{u}%
,\omega\right)  =H_{q}\left(  u,\omega\Phi_{u},\omega\right) \\
& =\omega\int\left(  1-q\Phi_{u}\right)  u^{2}dx\\
& =2\,\omega\,K_{q}\left(  u\right)  .
\end{align*}

Notice that, for $q=0,$ all these functionals reduce to the analogous ones for
the equation (\ref{KG}).

By the following proposition the study of the equations (\ref{a11}) and
(\ref{a222}) is reduced to an eigevalue problem of the type (\ref{JK}).

\begin{proposition}
\label{var}Let $q>0$ and $(u,\omega)\in H^{1}(\mathbb{R}^{3})\times\mathbb{R}$
be a solution of the eigenvalue problem%
\begin{equation}
J^{\prime}(u)=\omega^{2}K_{q}^{\prime}\left(  u\right)  .\label{ei}%
\end{equation}
Then $u,\varphi_{u},\omega$ solve (\ref{a11}) and (\ref{a222}).
\end{proposition}

Proof. First observe that $u,\varphi,\omega$ solve (\ref{a11}), (\ref{a222})
if and only if $\left(  u,\varphi\right)  $ is a critical point of the
functional
\begin{equation}
I_{\omega}(u,\varphi)=J(u)-F_{q}(u,\varphi,\omega)\label{ds}%
\end{equation}

namely if%
\begin{equation}
\frac{\partial I_{\omega}(u,\varphi)}{\partial u}=0,\text{ }\frac{\partial
I_{\omega}(u,\varphi)}{\partial\varphi}=0.\label{sat}%
\end{equation}

Now let $\left(  u,\omega\right)  $ be a solution of the eigenvalue problem
(\ref{ei}). Then clearly $u$ is a critical point of the functional
$u\rightarrow J(u)-\omega^{2}K_{q}\left(  u\right)  $ or equivalently, by
(\ref{sea}) and (\ref{ds}), a critical point of the functional%
\begin{equation}
u\rightarrow I_{\omega}(u,\varphi_{u})=J(u)-F_{q}(u,\varphi_{u},\omega
).\label{sata}%
\end{equation}

This means that%
\begin{equation}
\frac{\partial I_{\omega}(u,\varphi_{u})}{\partial u}+\frac{\partial
I_{\omega}(u,\varphi_{u})}{\partial\varphi}\varphi_{u}^{\prime}=0.\label{sate}%
\end{equation}
Since $\varphi_{u}$ solves (\ref{a222}), we have
\begin{equation}
\frac{\partial I_{\omega}(u,\varphi_{u})}{\partial\varphi}=0.\label{sati}%
\end{equation}
Then from (\ref{sate}) and (\ref{sati}) we get%
\begin{equation}
\frac{\partial I_{\omega}(u,\varphi_{u})}{\partial u}=0,\text{ }\frac{\partial
I_{\omega}(u,\varphi_{u})}{\partial\varphi}=0.\label{sato}%
\end{equation}
So by (\ref{sato}) we have that $u,\varphi_{u}$ solve (\ref{sat}).

$\square$

We shall show that if $q$ is small enough the eigenvalue problem (\ref{ei})
satisfies all the assumptions of the abstract theorem \ref{abstract}. More
precisely in this case we shall set%

\[
X=\left\{  u\in H^{1}(\mathbb{R}^{3}):u\ \text{is radially symmetric}\right\}
,
\]%
\[
\left\langle L_{1}u,u\right\rangle =\int\left(  \left\vert \nabla u\right\vert
^{2}+m^{2}u^{2}\right)  dx;\ N_{1}(u)=\int N(u)dx,
\]%
\begin{align*}
J(u)  & =\frac{1}{2}\left\langle L_{1}u,u\right\rangle +N_{1}(u)\\
& =\frac{1}{2}\int\left(  \left\vert \nabla u\right\vert ^{2}+m^{2}%
u^{2}\right)  dx+\int N(u)dx,
\end{align*}%
\begin{align*}
\left\langle L_{0}u,u\right\rangle  & =\int u^{2}dx,\\
\ K_{q}(u)  & =\frac{1}{2}\left\langle L_{0}u,u\right\rangle +N_{0}(u),\text{
}N_{0}(u)=-\frac{q}{2}\int\Phi_{u}u^{2}dx.
\end{align*}

\begin{lemma}
\label{pa}Assumptions (H1), (H2), (H3) are satisfied.
\end{lemma}

Proof. Arguing as in the proof of Lemma\ref{in} it can be proved that
assumption (H1) is satisfied and that $N_{1}^{\prime}$ satisfies (\ref{stare}).

Then, in order to complete the proof of (H2), we need to show that
$N_{0}^{\prime}$ is compact. First of all we look for a suitable expression
for $N_{0}^{\prime}.$

Observe that
\begin{equation}
K_{q}^{\prime}(u)=u+N_{0}^{\prime}(u).\label{la}%
\end{equation}

On the other hand by (\ref{bi}) and (\ref{sita})
\[
K_{q}\left(  u\right)  =F_{q}\left(  u,\Phi_{u},1\right)  .
\]
Then%
\begin{equation}
K_{q}^{\prime}\left(  u\right)  =\frac{\partial F_{q}\left(  u,\Phi
_{u},1\right)  }{\partial u}+\frac{\partial F_{q}\left(  u,\Phi_{u},1\right)
}{\partial\varphi}\Phi_{u}^{\prime}.\label{le}%
\end{equation}

Since $\Phi_{u}$ solves (\ref{b2}) and taking into account the definition
(\ref{be}) of $F_{q},$we have%
\begin{equation}
\frac{\partial F_{q}\left(  u,\Phi_{u},1\right)  }{\partial\varphi}=0,\text{
}\frac{\partial F_{q}\left(  u,\Phi_{u},1\right)  }{\partial u}=(1-q\Phi
_{u})^{2}u.\label{li}%
\end{equation}
So, comparing (\ref{le}), (\ref{li}), we have%
\begin{equation}
K_{q}^{\prime}\left(  u\right)  =(1-q\Phi_{u})^{2}u.\label{lo}%
\end{equation}
By (\ref{la}), (\ref{lo}) we get the following expression for $N_{0}^{\prime
}(u)$
\[
N_{0}^{\prime}(u)=(1-q\Phi_{u})^{2}u-u=q^{2}\Phi_{u}^{2}u-2q\Phi_{u}u.
\]
Then in order to show that $N_{0}^{\prime}$ is compact it is enough to prove
that the maps
\begin{equation}
u\rightarrow\Phi_{u}u\text{ and }u\rightarrow\Phi_{u}^{2}u\label{an}%
\end{equation}
are compact from $X$ to $X^{\prime}.$

Let
\[
u_{n}\rightharpoonup u_{0}\text{ weakly in }X.
\]
We shall prove first that $\Phi_{u_{n}}$ is bounded in $\mathcal{D}%
^{1,2}(\mathbb{R}^{3})$ and that, up to a subsequence,%
\begin{equation}
\Phi_{u_{n}}\rightharpoonup\Phi_{u_{0}}\text{ weakly in }\mathcal{D}%
^{1,2}(\mathbb{R}^{3}).\label{ui}%
\end{equation}
Since $\Phi_{u_{n}}$ solves
\begin{equation}
-\Delta\Phi_{u_{n}}+q^{2}u_{n}^{2}\Phi_{u_{n}}=qu_{n}^{2},\label{b3}%
\end{equation}

we have
\begin{equation}%
%TCIMACRO{\dint }%
%BeginExpansion
{\displaystyle\int}
%EndExpansion
\left\vert \nabla\Phi_{u_{n}}\right\vert ^{2}+q^{2}%
%TCIMACRO{\dint }%
%BeginExpansion
{\displaystyle\int}
%EndExpansion
\Phi_{u_{n}}^{2}u_{n}^{2}=q%
%TCIMACRO{\dint }%
%BeginExpansion
{\displaystyle\int}
%EndExpansion
\Phi_{u_{n}}u_{n}^{2}.\label{uo}%
\end{equation}
On the other hand%
\begin{equation}%
%TCIMACRO{\dint }%
%BeginExpansion
{\displaystyle\int}
%EndExpansion
\Phi_{u_{n}}u_{n}^{2}\leq\left\Vert \Phi_{u_{n}}\right\Vert _{L^{6}}\left\Vert
u_{n}\right\Vert _{L^{\frac{12}{5}}}^{2}.\label{ue}%
\end{equation}

Since $u_{n}$ is bounded in $X$, it is also bounded in $L^{\frac{12}{5}},$
then by (\ref{ue}) we have%
\begin{equation}%
%TCIMACRO{\dint }%
%BeginExpansion
{\displaystyle\int}
%EndExpansion
\Phi_{u_{n}}u_{n}^{2}\leq c_{1}\left\Vert \Phi_{u_{n}}\right\Vert _{L^{6}%
}.\label{ua}%
\end{equation}
From (\ref{uo}), (\ref{ua}) we easily get%
\[
\left\Vert \Phi_{u_{n}}\right\Vert _{\mathcal{D}^{1,2}}^{2}\leq c_{2}%
\left\Vert \Phi_{u_{n}}\right\Vert _{\mathcal{D}^{1,2}},
\]
from which we have that, up to a subsequence,%
\[
\Phi_{u_{n}}\rightharpoonup\Phi_{0}\text{ weakly in }\mathcal{D}%
^{1,2}(\mathbb{R}^{3})\text{.}%
\]
In order to prove (\ref{ui}) we have to show that $\Phi_{0}=\Phi_{u_{0}}$ i.e.
we show that $\Phi_{0}$ solves (\ref{b2}) with $u=u_{0}.$.

Let $v$ $\in C_{0}^{\infty}$ then, testing (\ref{b3}) on $v$ and passing to
the limit, we easily get%
\[
-\Delta\Phi_{0}+q^{2}u_{0}^{2}\Phi_{0}=qu_{0}^{2}.
\]
Then (\ref{ui}) is proved.

Now we prove that
\begin{equation}
u_{n}\Phi_{u_{n}}\rightarrow u_{0}\Phi_{u_{0}}\text{ in }L^{2}.\label{fa}%
\end{equation}
Let $\varepsilon,R>0$ and set%
\[
B_{R}=\left\{  x\in\mathbb{R}^{3}:\left\vert x\right\vert <R\right\}  ,\text{
}B_{R}^{c}=\mathbb{R}^{3}-B_{R}.
\]
Clearly we have%
\begin{equation}%
%TCIMACRO{\dint _{B_{R}^{c}}}%
%BeginExpansion
{\displaystyle\int_{B_{R}^{c}}}
%EndExpansion
\Phi_{u_{n}}^{2}u_{n}^{2}\leq\left(
%TCIMACRO{\dint _{B_{R}^{c}}}%
%BeginExpansion
{\displaystyle\int_{B_{R}^{c}}}
%EndExpansion
\left\vert u_{n}\right\vert ^{3}\right)  ^{\frac{2}{3}}\left(
%TCIMACRO{\dint _{B_{R}^{c}}}%
%BeginExpansion
{\displaystyle\int_{B_{R}^{c}}}
%EndExpansion
\Phi_{u_{n}}^{6}\right)  ^{\frac{1}{3}}.\label{cop}%
\end{equation}

Now we have (see \cite{Beres-Lions})%
\begin{equation}
\left\vert u_{n}(x)\right\vert \leq c_{1}\frac{\left\Vert u_{n}\right\Vert
_{H^{1}}}{\left\vert x\right\vert }\text{ in }B_{R}^{c}.\label{cap}%
\end{equation}
From (\ref{cop}) and (\ref{cap}) we get%
\begin{equation}%
%TCIMACRO{\dint _{B_{R}^{c}}}%
%BeginExpansion
{\displaystyle\int_{B_{R}^{c}}}
%EndExpansion
\Phi_{u_{n}}^{2}u_{n}^{2}\leq\left(  c_{1}\frac{\left\Vert u_{n}\right\Vert
_{H^{1}}}{R}\right)  ^{\frac{2}{3}}\left(
%TCIMACRO{\dint _{B_{R}^{c}}}%
%BeginExpansion
{\displaystyle\int_{B_{R}^{c}}}
%EndExpansion
\left\vert u_{n}\right\vert ^{2}\right)  ^{\frac{2}{3}}\left\Vert \Phi_{u_{n}%
}\right\Vert _{L^{6}}^{2}.\label{cip}%
\end{equation}

So, since $u_{n}$ is bounded in $H^{1}$ and $\Phi_{u_{n}}$ is bounded in
$\mathcal{D}^{1,2}(\mathbb{R}^{3})$ and hence in $L^{6},$ if we choose $R$
large enough, we get%
\begin{equation}%
%TCIMACRO{\dint _{B_{R}^{c}}}%
%BeginExpansion
{\displaystyle\int_{B_{R}^{c}}}
%EndExpansion
\Phi_{u_{n}}^{2}u_{n}^{2}<\varepsilon.\label{ciop}%
\end{equation}
Analogously, for $R$ large enough, we have%
\begin{equation}%
%TCIMACRO{\dint _{B_{R}^{c}}}%
%BeginExpansion
{\displaystyle\int_{B_{R}^{c}}}
%EndExpansion
\Phi_{u_{0}}^{2}u_{0}^{2}<\varepsilon\label{ciup}%
\end{equation}
and therefore%
\begin{equation}%
%TCIMACRO{\dint _{B_{R}^{c}}}%
%BeginExpansion
{\displaystyle\int_{B_{R}^{c}}}
%EndExpansion
\left\vert \Phi_{u_{n}}u_{n}-\Phi_{u_{0}}u_{0}\right\vert ^{2}<2\varepsilon
.\label{siup}%
\end{equation}
On the other hand%
\begin{align}%
%TCIMACRO{\dint _{B_{R}}}%
%BeginExpansion
{\displaystyle\int_{B_{R}}}
%EndExpansion
\left\vert \Phi_{u_{n}}u_{n}-\Phi_{u_{0}}u_{0}\right\vert ^{2}  & =%
%TCIMACRO{\dint _{B_{R}}}%
%BeginExpansion
{\displaystyle\int_{B_{R}}}
%EndExpansion
\left(  \Phi_{u_{n}}(u_{n}-u_{0})+u_{0}(\Phi_{u_{n}}-\Phi_{u_{0}})\right)
^{2}\leq\nonumber\\
& 2%
%TCIMACRO{\dint _{B_{R}}}%
%BeginExpansion
{\displaystyle\int_{B_{R}}}
%EndExpansion
\Phi_{u_{n}}^{2}(u_{n}-u_{0})^{2}+u_{0}^{2}(\Phi_{u_{n}}-\Phi_{u_{0}}%
)^{2}\nonumber\\
& \leq2\left\Vert \Phi_{u_{n}}\right\Vert _{L^{6}(B_{R})}^{2}\left\Vert
u_{n}-u_{0}\right\Vert _{L^{3}(B_{R})}^{2}\label{gua}\\
& +2\left\Vert u_{0}\right\Vert _{L^{6}(B_{R})}^{2}\left\Vert \Phi_{u_{n}%
}-\Phi_{u_{0}}\right\Vert _{L^{3}(B_{R})}^{2}.
\end{align}
The sequence $u_{n}$ weakly converges to $u_{0}$ in $H^{1}$, then it strongly
converges to $u_{0}$ in $L^{3}(B_{R}).$ So, since $\Phi_{u_{n}}$ is bounded in
$L^{6},$ we have%
\begin{equation}
\left\Vert \Phi_{u_{n}}\right\Vert _{L^{6}(B_{R})}\left\Vert u_{n}%
-u_{0}\right\Vert _{L^{3}(B_{R})}\rightarrow0.\label{pra}%
\end{equation}
On the other hand $\Phi_{u_{n}}\rightharpoonup$ $\Phi_{u_{0}}$ weakly in
$\mathcal{D}^{1,2}$ $\subset H_{loc}^{1}\subset\subset L_{loc}^{3},$ then we
have
\begin{equation}
\left\Vert \Phi_{u_{n}}-\Phi_{u_{0}}\right\Vert _{L^{3}(B_{R})}\rightarrow
0.\label{pre}%
\end{equation}

By (\ref{gua}), (\ref{pra}) and (\ref{pre}) we get%
\begin{equation}%
%TCIMACRO{\dint _{B_{R}}}%
%BeginExpansion
{\displaystyle\int_{B_{R}}}
%EndExpansion
\left\vert \Phi_{u_{n}}u_{n}-\Phi_{u_{0}}u_{0}\right\vert ^{2}\rightarrow
0.\label{siap}%
\end{equation}
Finally by (\ref{siup}) and (\ref{siap}) we get (\ref{fa}).

Following analogous arguments it can be shown that also the map $u\rightarrow
\Phi_{u}^{2}u$ is compact from $X$ to $X^{\prime}.$

Finally we prove that assumption (H3) is satisfied i. e. we prove that%
\[
\int\left(  1-q\Phi_{u}\right)  u^{2}dx\geq0.
\]

Arguing by contradiction assume that there is a region $\Omega$ where
$q\Phi_{u}>1$ and $q\Phi_{u}=1$ on $\partial\Omega.$ Clearly by (\ref{b2})
\[
-\Delta\left(  \Phi_{u}-\frac{1}{q}\right)  +q^{2}u^{2}\left(  \Phi_{u}%
-\frac{1}{q}\right)  =-\Delta\Phi_{u}+q^{2}u^{2}\Phi_{u}-qu^{2}=0.
\]
Then $v=\Phi_{u}-\frac{1}{q}$ solves the Dirichlet problem
\[
-\Delta v+q^{2}u^{2}v=0\text{ in }\Omega,\text{ }v=0\text{ on }\partial\Omega.
\]
Multiplying by $v$ and integrating in $\Omega$ we get
\[
\int_{\Omega}\left(  \left\vert \nabla v\right\vert ^{2}+q^{2}u^{2}%
v^{2}\right)  dx=0.
\]
Then $v=\Phi_{u}-\frac{1}{q}=0$ in $\Omega$ contradicting $q\Phi_{u}>1$ in
$\Omega$.

Finally observe that, if we take $u\neq0$ in all $\mathbb{R}^{3},$ then
\[
\int\left(  1-q\Phi_{u}\right)  u^{2}dx>0.
\]
In fact $\int\left(  1-q\Phi_{u}\right)  u^{2}dx=0$ would imply that $\Phi
_{u}=\frac{1}{q}$ a.e. in $\mathbb{R}^{3},$ contradicting $\Phi_{u}%
\in\mathcal{D}^{1,2}(\mathbb{R}^{3}).$

$\square$

\begin{lemma}
\label{bis} Assumption (\ref{HH}) is satisfied for $q$ sufficiently small.
\end{lemma}

Proof. Let $R>0$ and consider the map $u_{R}$ defined in (\ref{mapp}). As
shown in the proof of Lemma \ref{seguente}, we can choose $R$ be so large that%
\begin{equation}
\frac{J(u_{R})}{\frac{1}{2}%
%TCIMACRO{\dint }%
%BeginExpansion
{\displaystyle\int}
%EndExpansion
u_{R}^{2}}<m_{0}^{2}.\label{pri}%
\end{equation}
Now consider%
\begin{equation}
\frac{J(u_{R})}{K_{q}(u_{R})}=\frac{J(u_{R})}{\frac{1}{2}%
%TCIMACRO{\dint }%
%BeginExpansion
{\displaystyle\int}
%EndExpansion
u_{R}^{2}-\frac{q}{2}\int\Phi_{u_{R}}u_{R}^{2}}.\label{prix}%
\end{equation}
So, by (\ref{prix}), we get that assumption (\ref{HH}) is satisfied if we show
that%
\begin{equation}
q\int\Phi_{u_{R}}u_{R}^{2}\rightarrow0\text{ for }q\rightarrow0.\label{fina}%
\end{equation}
Since $\Phi_{u_{R}}$ depends on $q$ a little work is needed to prove
(\ref{fina}).

Since $\Phi_{u_{R}}$ solves (\ref{b2}) with $u=u_{R},$ we have%

\[
\left\Vert \Phi_{u_{R}}\right\Vert _{\mathcal{D}^{1,2}}^{2}+q^{2}\int
u_{R}^{2}\Phi_{u_{R}}^{2}=q\int u_{R}^{2}\Phi_{u_{R}}\leq
\]%
\begin{equation}
\leq q\left\Vert u_{R}\right\Vert _{L^{\frac{12}{5}}}^{2}\left\Vert
\Phi_{u_{R}}\right\Vert _{L^{6}}\label{again}%
\end{equation}
and then
\[
\frac{\left\Vert \Phi_{u_{R}}\right\Vert _{\mathcal{D}^{1,2}}^{2}}{\left\Vert
\Phi_{u_{R}}\right\Vert _{L^{6}}}\leq q\left\Vert u_{R}\right\Vert
_{L^{\frac{12}{5}}}^{2}.
\]
Then, since $\mathcal{D}^{1,2}$ is continuously embedded into $L^{6}$, we
easily get
\begin{equation}
\left\Vert \Phi_{u_{R}}\right\Vert _{\mathcal{D}^{1,2}}\leq cq\left\Vert
u_{R}\right\Vert _{L^{\frac{12}{5}}}^{2},
\end{equation}
where $c$ is a positive constant. Then, using again (\ref{again}), we get
\[
q\int u_{R}^{2}\Phi_{u_{R}}\leq q\left\Vert u_{R}\right\Vert _{L^{\frac{12}%
{5}}}^{2}\left\Vert \Phi_{u_{R}}\right\Vert _{L^{6}}\leq cq^{2}\left\Vert
u_{R}\right\Vert _{L^{\frac{12}{5}}}^{4}.
\]
From which we get (\ref{fina}).

$\square$

Finally we are ready to conclude the proof of Theorem
\ref{main-theorem copy(1)}.

Proof of Theorem \ref{main-theorem copy(1)}.

By Lemma \ref{pa} the assumptions (H1), (H2), (H3) of the Theorem
\ref{abstract} are satisfied. Moreover by Lemma \ref{bis} there exists
$q^{\ast}>0$ such that for $0<q<q^{\ast}$ also assumption (\ref{HH}) is
satisfied. Then we can use Theorem \ref{final} and we get that there exists
$q^{\ast}>0 $ such that for $0<q<q^{\ast}$ there exists a non empty, open
subset $\Sigma_{q}\subset\mathbb{R}$ such that for any $\sigma\in\Sigma_{q}$
problem (\ref{ei}) has a solution $(u,\omega)$ with charge $H_{q}%
(u,\omega)=\sigma.$ Moreover such a solution minimizes the energy $\tilde
{E}_{q}\left(  u,\omega\right)  $ on the states $\left(  u,\omega\right)  $
having charge $H_{q}(u,\omega)=\sigma.$ Then, by Proposition \ref{var},
$u,\omega,$ $\varphi_{u}=\omega\Phi_{u}$ solve (\ref{a11}), (\ref{a222}).

$\square$

\section{Appendix.}

Let assumption Wiii) (b) be satisfied i.e. we assume that there exists
$s_{1}>s_{0}$ such that $N^{\prime}(s_{1})\geq0.$

Set%
\begin{equation}
\tilde{N}(s)=\left\{
\begin{array}
[c]{l}%
N(s)\text{ for }s\leq s_{1}\\
N^{\prime}(s_{1})s+c_{1}\text{ for }s\geq s_{1}%
\end{array}
\right. \label{immm}%
\end{equation}

where
\[
c_{1}=N(s_{1})-N^{\prime}(s_{1})s_{1}%
\]

Set
\begin{equation}
\tilde{W}(s)=\frac{m^{2}}{2}s^{2}+\tilde{N}(s)\label{imm}%
\end{equation}

By the following proposition we can replace in (\ref{static}) $W^{\prime}(s) $
with $\tilde{W}^{\prime}(s)$

\begin{proposition}
\label{max}Let $m^{2}\geq\omega^{2}.$ Then for any solution $u\in H^{1}$ of
the equation
\begin{equation}
-\Delta u+\tilde{W}^{\prime}(u)=\omega^{2}u\label{one}%
\end{equation}
we have
\[
u\leq s_{1}%
\]

\end{proposition}

Proof. Let $u\in H^{1}$ be a solution of (\ref{one}) and set
\[
u=s_{1}+v.
\]
We want to show that $v\leq0.$ Arguing by contradiction, assume that
\[
\Omega=\left\{  x:v(x)>0\right\}  \neq\varnothing.
\]
Then, multiplying both members of (\ref{one}) by $v$ and integrating on
$\Omega,$we have%
\[%
%TCIMACRO{\dint \limits_{\Omega}}%
%BeginExpansion
{\displaystyle\int\limits_{\Omega}}
%EndExpansion
\left\vert \nabla v\right\vert ^{2}+\tilde{W}^{\prime}(s_{1}+v)v-\omega
^{2}(s_{1}+v)v=0.
\]
So, using (\ref{imm}), we have%
\[%
%TCIMACRO{\dint \limits_{\Omega}}%
%BeginExpansion
{\displaystyle\int\limits_{\Omega}}
%EndExpansion
\left\vert \nabla v\right\vert ^{2}+\tilde{N}^{\prime}(s_{1}+v)v+(m^{2}%
-\omega^{2})(s_{1}+v)v=0
\]
which, by (\ref{immm}), becomes%
\begin{equation}%
%TCIMACRO{\dint \limits_{\Omega}}%
%BeginExpansion
{\displaystyle\int\limits_{\Omega}}
%EndExpansion
\left\vert \nabla v\right\vert ^{2}+N^{\prime}(s_{1})v+(m^{2}-\omega
^{2})(s_{1}+v)v=0.\label{exp}%
\end{equation}

Since
\[
N^{\prime}(s_{1})\geq0\text{ and }m^{2}\geq\omega^{2},
\]
expression (\ref{exp}) gives%
\[
v=0\text{ in }\Omega,
\]
contradicting the definition of $\Omega.$

$\square$

\bigskip

\end{document}